# The Value of Recycling for Low-Carbon Energy Systems - a Case Study of Germany's Energy Transition

Felix Kullmann[1,*], Peter Markewitz[1], Leander Kotzur[1], Detlef Stolten[1,2]

[1] Institute of Energy and Climate Research, Techno-economic Systems Analysis (IEK-3) Forschungszentrum Jülich GmbH, Wilhelm-Johnen-Str., 52428 Jülich, Germany

[2] Chair of fuel cells, RWTH Aachen University, c/o Institute of Energy and Climate Research, Techno-economic Systems Analysis (IEK-3) Forschungszentrum Jülich GmbH, Wilhelm-Johnen-Str., 52428 Jülich, Germany

[*] Corresponding author. E-mail address: f.kullmann@fz-juelich.de (F. Kullmann).

## Abstract

To achieve climate neutrality, synergies between circular economy and reduction of greenhouse gas emissions must be strengthened. Previously idle emission reduction potentials of resource efficiency are to be exploited. Since all potentially possible emission reduction measures are linked by interactions, the evaluation of a single measure in terms of cost efficiency, effectiveness, and compliance with climate protection targets is very complex and requires a model-based analysis that takes the entire energy system into account.

This work advances an energy system model for Germany so that through comprehensive modeling of industrial processes and implementation of recycling options, the impact of recycling measures in the context of national greenhouse gas mitigation strategies can be analyzed.

The scenario evaluation shows that different recycling strategies have large effects on the German energy system. Without recycling energy demand in 2050 will increase by more than 300 TWh and cost of transformation will rise by 85% compared to a reference scenario. On the other hand, if maximum recycling rates can be achieved, costs of transformation can be reduced by 26% until 2050. Recycling is an essential and cost-efficient greenhouse gas reduction strategy for future low-carbon energy system designs.

## 1 Introduction

With the "European Green Deal," the European Commission presents its roadmap for transformation for a sustainable future EU economy [1]. The European Commission equates efforts that lead from a linear to a circular industry with efforts for more climate protection. With the resulting "new action plan for the circular economy," which was published in 2020, the European Commission undertakes to analyze "[...] the impact of the circular economy on climate change mitigation [...]" [2, p. 19]. "To achieve climate neutrality, the synergies between the circular economy and the reduction of greenhouse gas (GHG) emissions must be strengthened. " [2, p. 19]. As early as 2016, the federal government of Germany set with the "Climate Protection Plan 2050" [3] long-term targets for the reduction of greenhouse gas emissions by 2050, based on the "Paris Agreement" [4]. Concrete measures for achieving the targets were presented in 2020 in the "Climate Protection Program 2030" of the Federal Government for the Implementation of the "Climate Protection Plan 2050" [5]. It states "[...] that the targets will be achieved most cost-effectively if they can be realized across all sectors." [5, p. 8]. In particular, the principle of circular economy is prominently represented in the measures presented there. With the help of the circular economy, previously unexploited emission reduction potentials of resource efficiency are to be tapped. The German Resource Efficiency Program III of 2020 also emphasizes the link between the circular economy and climate protection. "[I]t will not be possible to meet the target [...] set out in the Paris Agreement on climate protection [...]" without raw material efficiency measures [...]." [6, p. 6]

All potentially possible measures to reduce emissions are linked by interactions between the individual sectors in the energy system. Due to the resulting interactions, the assessment in terms of cost efficiency, effectiveness and compliance with climate protection targets of a single measure is not possible in most cases.

Due to the high degree of complexity, a comprehensive assessment that takes all interactions into account is only possible with a model-based analysis. For this purpose, the use of energy system models that depict the entire energy supply across all sectors is appropriate. In the past, fields of action such as recycling and material recycling of $CO_2$ were understood less as energy topics and more as resource topics. Consequently, they are only rudimentarily represented in existing energy system models, or not at all. Therefore, the focus of the present work is the extension of an existing national energy system model, which allows to estimate the contribution of recycling with regard to the reduction of climate gas emissions.

*Motivation and state of the art*

Energy system models are used to analyze future energy supply and demand structures and to assess impacts of policies on an energy system [7]. Although national energy systems will have to change drastically in order to guarantee an emission-free energy supply in the future, the concrete design of a future energy system is not yet clear [8]. There are a number of studies that use energy system models to create national scenarios for Germany, showing possible transformation paths towards an energy system in 2050. In most scenarios, recycling, in the sense of resource efficiency or recycling of $CO_2$ streams for material use, is not considered at all or is insufficiently considered. This means, that an assessment of the extent to which these measures are part of efficient greenhouse gas reduction strategies and what effects their absence implies for the overall German energy system has not yet been made. Furthermore, it is not clear, whether recycling measures are a cost-efficient option within an overall mitigation strategy.

The International Energy Agency (IEA) includes material efficiency and recycling strategies in its global future scenarios that lead to reductions in $CO_2$-emissions and energy use [9]. The study "The Circular Economy - A Powerful Force for Climate Mitigation" [10] concludes that demand-side measures can reduce $CO_2$-emissions from the European industrial sector by almost 300 million tons per year by 2050. Recycling opportunities alone contribute to about 60% of these savings. At the national level, the Renewable Energies Research Association (FVEE) emphasizes that recycling processes and economical use of materials are prerequisites for realizing a low-emission energy system in Germany [11]. Gerbert et al. [12] conclude that higher recycling rates of non-ferrous metals in Germany could lead to savings of up to 2 Mt $CO_2$-eq/a. Another study [13] puts the GHG emission savings from scrap-based nonferrous metal production in Germany at 7 Mt, with an increasing potential by 2050. To assess this energy savings and $CO_2$ reduction potential, these studies use exogenously predefined recycling rates. This practice provides little insight into whether recycling is the cost-optimal choice as a GHG reduction strategy. In addition, this approach does not allow the effects in interaction with other $CO_2$ reduction measures in the overall energy system to be illuminated.

Against this background, the present work aims to further develop an existing energy system model in such a way that, through comprehensive modeling of industrial processes and the implementation of recycling options, consistent national greenhouse gas reduction strategies for Germany can be investigated and, in addition, the effects of these measures on the overall energy system can be quantified and evaluated. The focus of the implementation is on potential recycling processes. The following research questions will be answered in this paper.

1. *What does a cost-optimal transformation of the German industrial sector look like in the context of an overall national greenhouse gas mitigation strategy?*
2. *What is the value of recycling for achieving the national greenhouse gas reduction targets for Germany?*

## 1.1 Circular economy of the industry in current energy scenarios

Since the aim of this work is to evaluate an integrated analysis of the effects of recycling and material use in the German energy system, the scenario selection in the following chapter is limited to national energy scenarios for Germany. Further criteria are that the time horizon of the selected scenarios should include the year 2050, and that the objective of the scenarios should be focused on the implementation of the climate protection goals of the German government [3,5]. The core results of the individual scenarios are presented in order to be able to classify the results of this work afterwards. A total of 30 scenarios from the last ten years are analyzed.

*Trends and developments in the scenarios studied*

A summary of the greenhouse gas reduction scenarios examined can be found in Table 1-1 which shows that the older scenarios in particular do not take recycling into account. In addition, the industrial sector is usually insufficiently represented, so that in some cases no statement can be made about the distribution of future final energy consumption on specific processes.

**Table 1-1 German energy system scenarios examined regarding implementation of recycling and industry (-: not considered; x: crude implementation; xx: detailed implementation).**

| ID | Year | Source | German energy system scenarios | Recycling | Industry |
|---|---|---|---|---|---|
| A | 2009 | [14] | Long-term scenarios (2009) | - | x |
| B1-2 | 2009 | [15] | Model Germany | - | x |
| C | 2010 | [16] | Energy target 2050: 100% from renewable energies | - | x |
| D1 | 2012 | [17] | Long-term scenarios (2012) | - | x |
| E | 2014 | [18] | Development of the energy markets - Energy reference forecast | x | x |
| F | 2014 | [19] | Energy transition business model | - | - |
| G | 2015 | [20] | Climate protection scenario 2050 | x | x |
| H | 2016 | [21] | Energy system 2050 | - | x |
| I | 2016 | [22] | The energy transition after COP 21 - Current scenarios for German energy supply | - | - |
| J | 2016 | [23] | Sector coupling through the energy transition | - | - |
| K1-2 | 2017 | [24] | Long-term scenarios (2017) | x | xx |
| L | 2017 | [25] | Shaping the path to a greenhouse gas-neutral Germany in a resource-conserving way (2017) | xx | xx |
| M1-3 | 2018 | [12] | Climate paths for Germany | - | xx |
| N1 | 2018 | [26] | Cost-efficient sector coupling | x | xx |
| O1-5 | 2018 | [27] | dena - Lead study on integrated energy transition | x | xx |
| P | 2019 | [28] | Shaping the path to a greenhouse gas-neutral Germany in a resource-conserving way (2019) | xx | xx |
| Q2 | 2019 | [29] | Paths for the energy transition | x | x |
| R | 2020 | [30] | Paths to a climate-neutral energy system | - | x |
| S | 2020 | [31] | Hydrogen Roadmap North Rhine-Westphalia | x | xx |
| T | 2020 | [32] | Climate-neutral Germany | x | xx |

In the more recent energy scenarios, the level of detail of industry mapping is increasing. However, in most cases, simulation models are used that update the final energy demand of industry and use it as an input parameter for optimizing the supply of energy sources. Thus, no statement can be made about the cost optimality of the overall energy system. No information about the cost efficiency and effectiveness for reducing greenhouse gas emissions of recycling measures can be retrieved across all scenarios examined. This is partly due to the fact that in the scenarios in which recycling is considered, only specific recycling rates are exogenously specified and are not part of an endogenous cost optimization.

For comparison purposes, Figure 1-1 shows the scenarios examined and the final energy consumption of the industrial sector reported in each case. The scenarios are sorted according to their year of publication and conclude on the right-hand side with the status quo of the year 2019 [14]. It is striking that the older scenarios show lower energy consumption on average than the newer scenarios. One reason for this is that in these scenarios, energy efficiency in the industrial sector was assessed as one of the most important and economic measures. Innovative technologies or a use of novel energy sources was not considered at all in many older scenarios. This is also confirmed by the fact that hydrogen plays a role as an energy carrier in industry in almost all recently published scenarios.

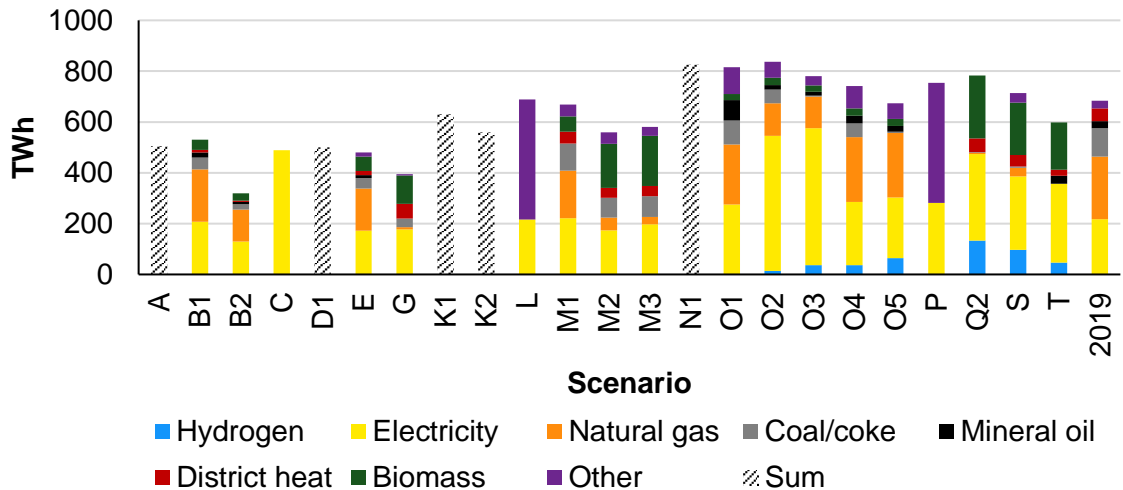

**Figure 1-1 Industrial final energy demand in selected scenarios in 2050 (2019 values from [14,15])**

It can also be observed that older scenarios often did not include a detailed statement about the final energy consumption in the industrial sector. This is because many of the models that served as the basis for these scenarios did not include a detailed implementation of the industrial sector and thus no statement could be made about the energy carrier mix in industry. In addition, a consequence of this is that no analysis of recycling measures can be made without a detailed mapping of industry. This observation fits with the fact that recycling is only accounted for in more recent national energy scenarios.

## 1.2 Energy system and material flow models

Based on 55 material flow models and 5 energy system models, a previous publication first discusses the individual model worlds separately and explains their characteristics, and then compares both model classes [16]. The problem of implementing material flows in energy system models, for the representation of recycling measures, is summarized again here. The following points can be made:

- The energy world is inadequately represented in material flow models and, conversely, material flows are inadequately represented in energy system models.
- Recycling measures cannot currently be adequately represented in energy system models.
- In many cases, the industrial sector is insufficiently represented. Strategies relating to recycling cannot be derived.
- Cost efficiency plays no role in almost any scenario.

- There is no consistent and cost-effective overall assessment of recycling in the context of the entire energy system. As a rule, the existing studies are merely simulations without analysis of the costs.
- Recycling is only taken into account, if at all, through exogenously set assumptions. In most scenarios, the topic of recycling is not addressed at all.

# 2 Methodology

## 2.1 Energy system model NESTOR

NESTOR stands for National Energy System Model with Sector Coupling and refers to an optimization model that maps the German energy supply system across all sectors. The fundamentals of this model are taken from Lopion [17] and Kotzur [18] and are briefly presented in supplementary material A. NESTOR is an integrated national energy system model (including all energy supply and demand sectors), which optimizes the German energy system design of a future target year (e.g., 2050) based on cost-efficiency. Additionally, the model optimizes the transformation path until the selected target year using a myopic approach.

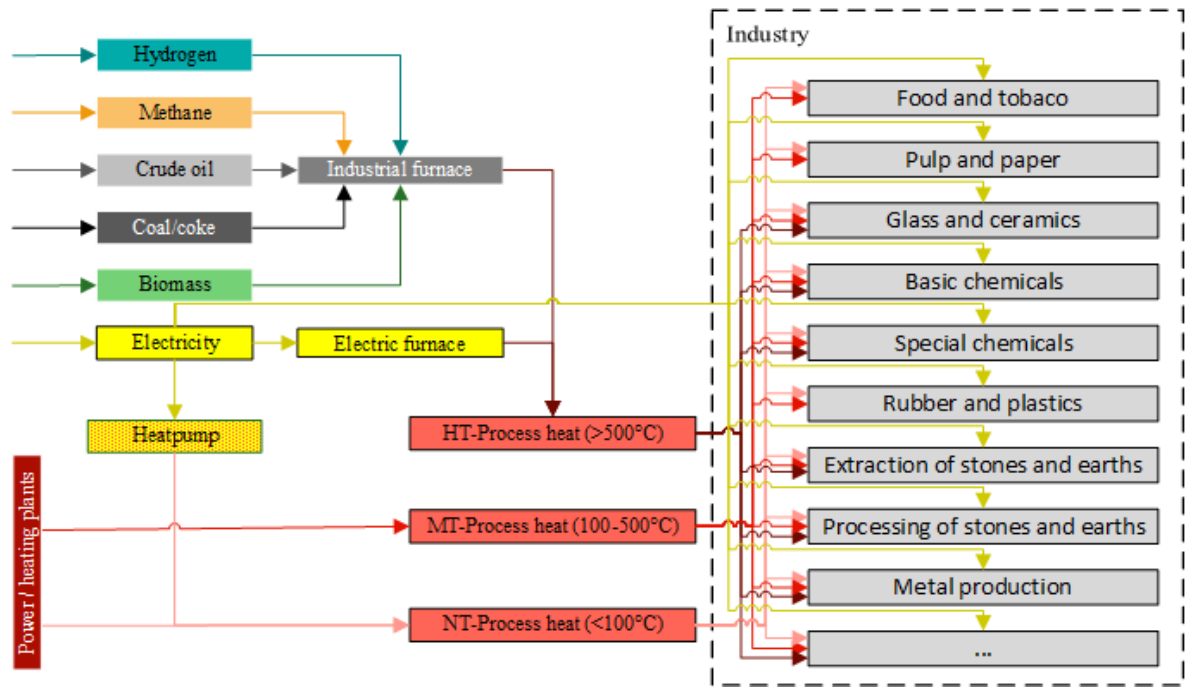

**Figure 2-1 Overview of how the industry sector is embedded in the energy system model**

The different industries are embedded in NESTOR as shown in Figure 2-1. Process heat and electricity for the processes are produced within the energy system. Which energy carrier is used is a result of the cost-optimization. For metal and non-ferrous metal production, basic chemicals, pulp and paper, glass production and processing of stones and earths, production processes have been modelled on a more detailed level. This enables to also include recycling technologies. Figure 2-2 shows how the methanol production is embedded in NESTOR as an example of the detailed process implementation.

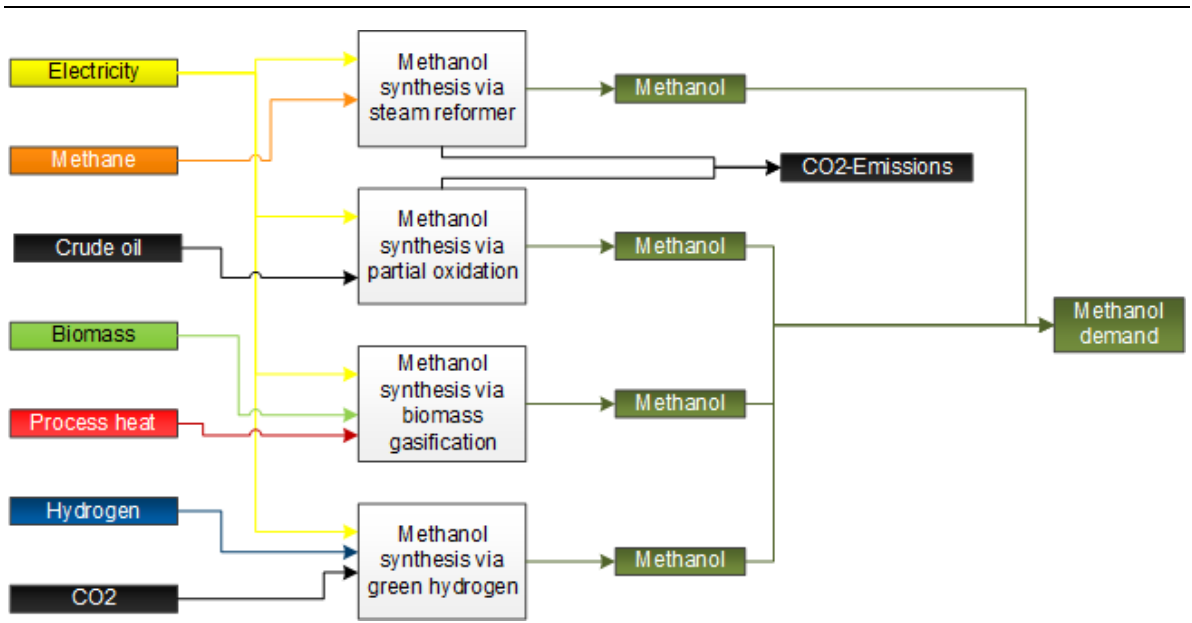

**Figure 2-2 Overview of how methanol production is embedded in the energy system model**

Several methanol production processes are implemented in NESTOR, which all can produce methanol and satisfy the demand. All energy carriers are generated endogenously and provide the link to the energy system. The corresponding parameters used to implement the industrial processes on a detailed level can be found in supplementary material B.

## 2.2 Recycling

With the model-based implementation of detailed process illustrations a fundamental prerequisite for the analysis of recycling is fulfilled. Another important component is the estimation of future available secondary raw material quantities in order to be able to evaluate the possible recycling potential. In the following, a procedure is presented that enables the estimation of secondary raw material quantities, which in turn are set as input parameters in the NESTOR model.

### 2.2.1 Estimation of material flows available in the future

For a detailed consideration of the greenhouse gas reduction potentials of recycling, it is not only necessary to consider current and future recycling processes and their techno-economic parameters, but also to be able to quantify the material flows that will provide raw materials for recycling in the future. For this estimation, the methods of material flow analysis are used in this work (cf. [19–23]). These models are usually based on simple balance sheet calculations. This approach is used to make a statement about what proportion of the waste/scrap quantities produced today will be available to the energy system in future. For the analysis of future waste/scrap flows, it is particularly important to know what proportion will remain in the anthropogenic stock for how long until it flows back into the overall system and is available for recycling. For this purpose, most material flow models use frequency distributions. The basic idea is that depending on which material is involved and in which anthropogenic area it is used, the residence time in this area changes. In order to be able to quantify future material flows with the energy system model, it is necessary to take these frequency distributions into account. This procedure is described in more detail below.

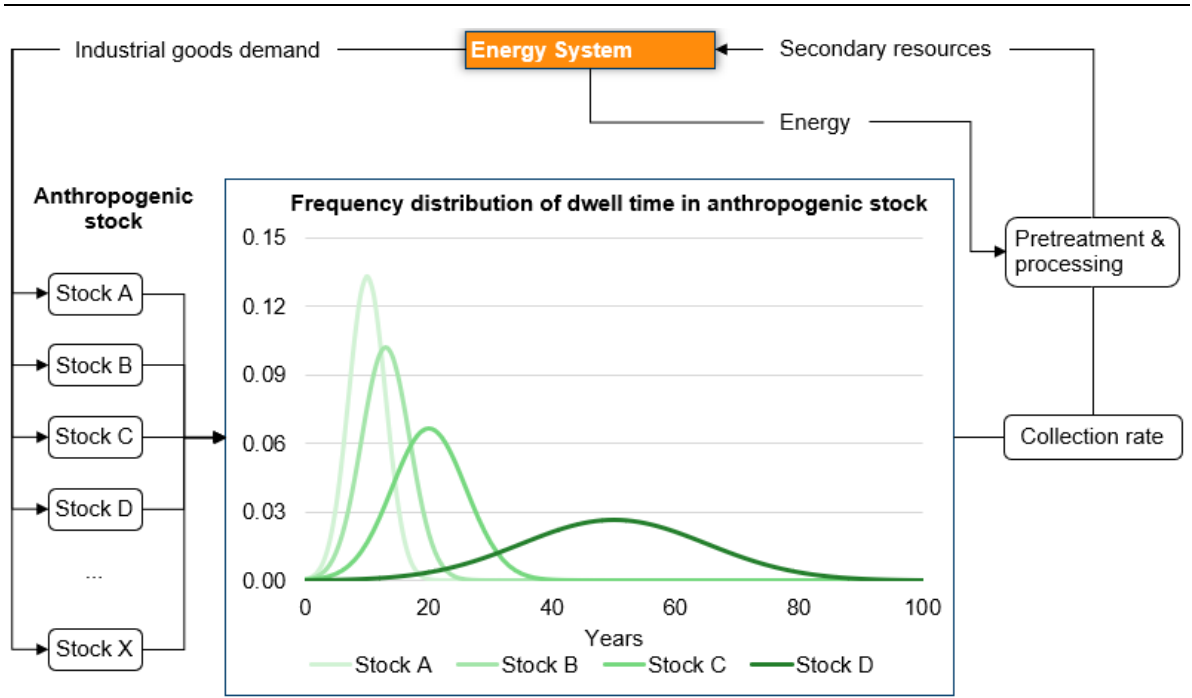

**Figure 2-3 Schematic representation of anthropogenic stocks in the energy system**

As shown in Figure 2-3, the respective demand for goods of the energy system is first divided into different anthropogenic stocks. Within these stocks, the amount of material bound in products remains for a specific time. For example, materials in stock A are sent to waste recycling after a relatively short lifetime of 10 years (e.g. plastic consumer goods, etc.), whereas materials in stock D have a much longer residence time and are only returned to waste recycling after an average of 50 years (e.g. construction industry, infrastructures, etc.). Following this approach the amount of waste of a material available to the energy system each year can be estimated (cf. Eq. 2-1). The amount of secondary raw material $b_{x,k,t}$ that will be produced in year $x$ from the amount of $D_{t,k}$ that was produced in year $t$ and goes to the anthropogenic stock $k$ can be calculated assuming the average retention time $\mu$ and the standard deviation $\sigma$.

$$b_{x,k,t} = \left( \frac{1}{\sigma_k \sqrt{2\pi}} e^{\frac{1}{2}\left(\frac{(x-t)-\mu_k}{\sigma_k}\right)^2} \right) * D_{t,k} \qquad \text{Eq. 2-12}$$

For example, if 10 Mt of an industrial good has entered the construction sector in 2020, then assuming a normal distribution with an assumed mean residence time of 50 years and a standard deviation of 15 years, there will be 7.6 kt of scrap in 2030 and 109.3 kt of scrap in 2050. Scrap flows can be quantified for all years following the point in time when a certain amount of this industrial good is produced and enters the individual anthropogenic stock. Thus, in an optimization year, the cumulative sum of the inputs of a material into the respective anthropogenic stock in all preceding years results in the theoretically available secondary raw material quantity (cf. Eq. 2-2).

$$B_{x,k} = \sum_{i=t}^{x-1} \left( \left( \frac{1}{\sigma_k \sqrt{2\pi}} e^{\frac{1}{2}\left(\frac{(x-i)-\mu_k}{\sigma_k}\right)^2} \right) * D_{i,k} \right) \qquad \text{Eq. 2-2}$$

For the years within the transformation path analysis (2020-2050), the quantities result endogenously from the model. However, these quantities must be supplemented in each case by the quantity of scrap that arose from the production of a good before the optimized transformation path and only reached its maximum residence time in the anthropogenic stock during the years of the transformation path analysis. Since these estimates are always the theoretical maximum amount of scrap, values for the yield and collection rates are additionally assumed for each substance, depending on the anthropogenic stock in which it resided. These assumptions can be found in supplementary material C.

**2.2.2 Endogenous recycling rates as part of model optimization.**

The recycling rates for a given material in NESTOR are not exogenously specified but are themselves part of the model solution and therefore subject to cost optimality. Costs, as well as material and energy demands, are applied to the entire recycling chain and for the recycling of a material. Thus, the processes for the primary production of a material compete with the processes for the reuse of scrap in the cost-optimal system. Figure 2-4 illustrates this process.

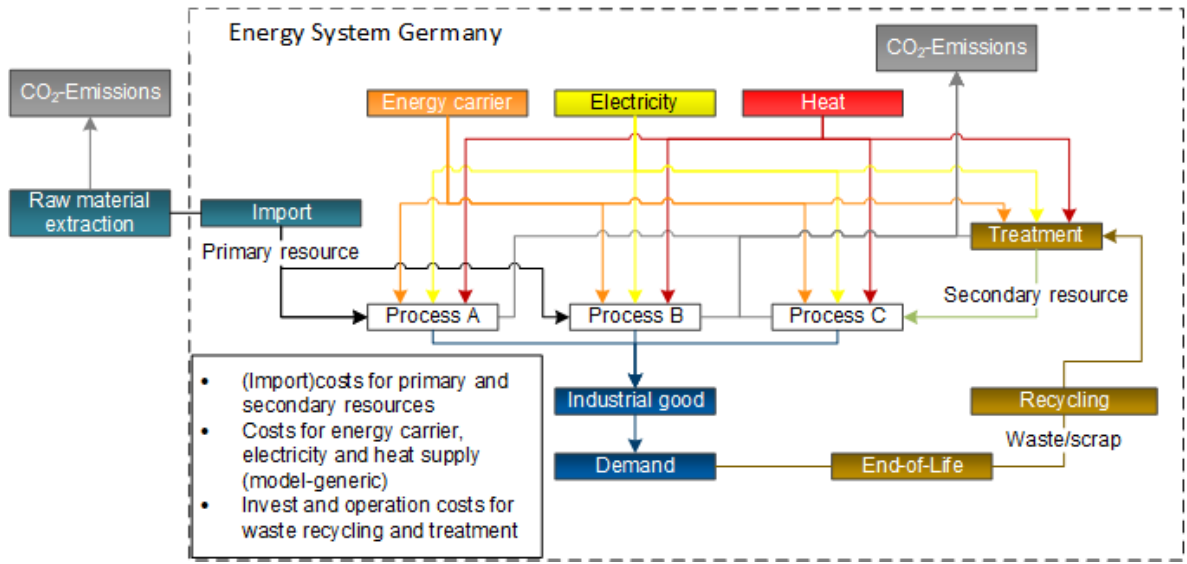

**Figure 2-4 Schematic representation of the raw material paths in the energy system model**

The energy system model is given an exogenous demand for goods, which must be met by the industrial sector represented. The processes by which the industrial goods are to be produced are fixed and part of the optimization result. Accordingly, by the year 2050, there are no exogenously specified proportions for the processes that can produce a particular industrial good. As a model result, the demand for industrial goods can be met entirely by a single process as well as by any ratio of the different processes. The processes are not only in cost competition with each other, but with all technologies of the whole energy system in order to obtain an economic optimum. For this purpose, it is necessary that both costs and specific $CO_2$ emissions are applied for the different process paths. This includes the direct costs (investment and operating costs) of the process for producing an industrial good as well as the costs of the domestic upstream chains. These include import costs for primary raw materials, procurement costs for electricity, heat and other energy sources. In the case of imports of the energy carriers required for the process, these are exogenously specified import costs (e.g. hydrogen import costs). Other energy carriers converted within the national boundary of Germany (e.g. electricity and process heat), are freely optimized in the energy system model. This means that the electricity prices result endogenously at each time step and are part of the optimization. In addition, no specific recycling rates are exogenously set. Due to the competition with processes of primary raw material processing, recycling quotas result based on the share of a recycling process in the cost-optimal solution. For recycling processes, the costs of the required secondary raw materials are of particular importance, but these must also be exogenously specified.

# 3 Scenarios

The development of the industrial sector takes place in the context of the overall energy system. Due to the systemic interactions in the energy system model, the optimization results in the industrial sector cannot be evaluated in isolation from the developments in other sectors. Therefore, an overview of the general assumptions concerning the overall system and the individual sectors, which underlie all scenarios, is given first.

## 3.1 Scenario description

As a starting point for the two recycling scenarios studies in this paper, a baseline scenario is created (*REF95*). In *REF95,* $CO_2$ emissions must be reduced by 95% in 2050 compared to 1990 (see Appendix A for basic assumptions of the energy system model). Apart from the phase-out of coal and nuclear power generation (according to [24,25]), no other limitation is placed on the model. Recycling rates are fixed at today's levels in the individual industrial sectors until 2050. Novel recycling processes (e.g., chemical recycling of plastic waste) are available from 2040. Otherwise, no further restrictions are applied. The base scenario *REF95* is used for comparison with the other scenarios and thus serves for evaluation and classification. A comparative overview with the most important assumptions of the scenarios can be found in Table 3-1.

**Table 3-1 Comparison of selected assumptions for the scenarios created**

|  | *REF95* | *w/o Rec* | *RecX* |
|---|---|---|---|
| **GHG reduction compared to 1990** | -95% | -95% | -95% |
| **Recycling rate in 2050** | fixed to 2020 | - | maximum |
| **New rec. technologies** | as of 2040 | - | as of 2035 |
| **Imported Hydrogen €/kWh** | 0.10 | 0.10 | 0.10 |
| **Imported Power-to-Liquid €/kWh** | 0.16 | 0.16 | 0.16 |
| **Steel scrap €/kg** | 0.24 | - | 0.24 |
| **Aluminum scrap €/kg** | 0.35 | - | 0.35 |
| **Copper scrap €/kg** | 2.65 | - | 2.65 |
| **Zinc scrap €/kg** | 1.47 | - | 1.47 |
| **Waste glass €/kg** | 0.05 | - | 0.05 |
| **Waste paper €/kg** | 0.05 | - | 0.05 |

The scenario *w/oRec* represents an extreme case in which it is assumed that no recycling is possible until the year 2050. Thus, the greenhouse gas reduction must take place completely without recycling measures. Compared to the reference scenario (*REF95*), this scenario, which is similar to a value-off analysis, can be used to analyze the value of recycling measures. For the scenario *RecX*, the model is given the opportunity to fully exploit the theoretical maximum recycling rates.

## 3.2 Scenario results

The following is a comparison of the recycling scenarios. All scenarios are mirrored against the *REF95* reference scenario. The focus is on the effects that become apparent in the overall system and, in particular, in the industrial sector. Primary energy consumption, electricity consumption and generation, and hydrogen demand and generation are used as indicators for overall system effects. The comparative analysis of the industrial sector is based on final energy consumption, hydrogen demand, and selected key technologies. Subsequently, the system costs of the scenarios are mirrored on the reference scenario *REF95 in order* to be able to make statements about the cost efficiency of certain measures or strategies.

### 3.2.1 Impact of recycling on the total energy system design

Figure 3-1 shows the development of primary energy consumption up to the year 2050 for the reference scenario. Due to the restrictive $CO_2$ limits, fossil energy sources are completely substituted by renewable energy sources (only fossil-based non-energetic demand in 2050). A change in recycling rates has a major impact on the total demand. Without recycling, the absolute primary energy consumption increases by approx. 300 TWh in 2050, which is an increase of 14% compared to *REF95*. On the other hand, increased recycling leads to an approx. 280 TWh lower primary energy demand (-12% compared to *REF95*).

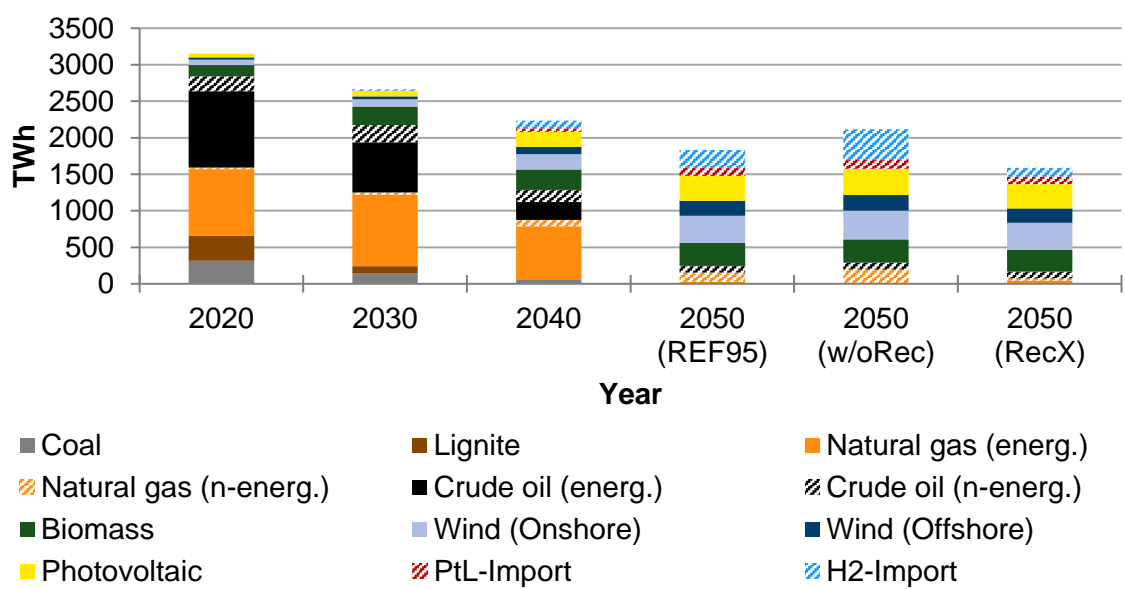

**Figure 3-1 Development of primary energy demand in reference scenario, with 2050-values across all scenarios**

Clear differences can be observed in all scenarios with regard to the necessary hydrogen import; higher recycling rates minimize dependency from energy imports.

The scenarios also differ in terms of electricity demand (see Figure 3-2). Whereas in the first few years low-cost energy efficiency potentials in the buildings sector will be exploited and result in lower electricity demand, the demand for electricity will double by 2050.

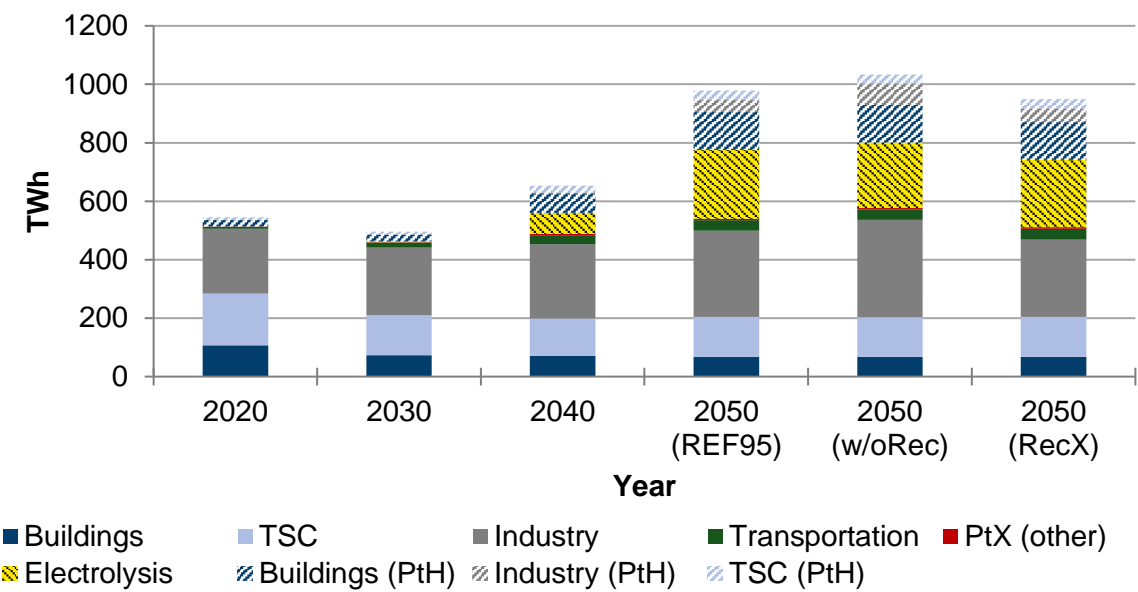

**Figure 3-2 Development of electricity demand in reference scenario, with 2050-values across all scenarios (TSC: Trade, Service, Commerce)**

This is mainly explained by the increased use of Power-to-X technologies. In the reference scenario, more than 236 TWh of electricity is required for domestic hydrogen production with electrolysers. Another 43 TWh of electricity is used to provide process heat (PtH) in industry. Without recycling measures, the electricity demand of industry increases by an additional 55 TWh in 2050.

This has implications for both installed capacity and electricity generation in 2050 (see Figure 3-3).

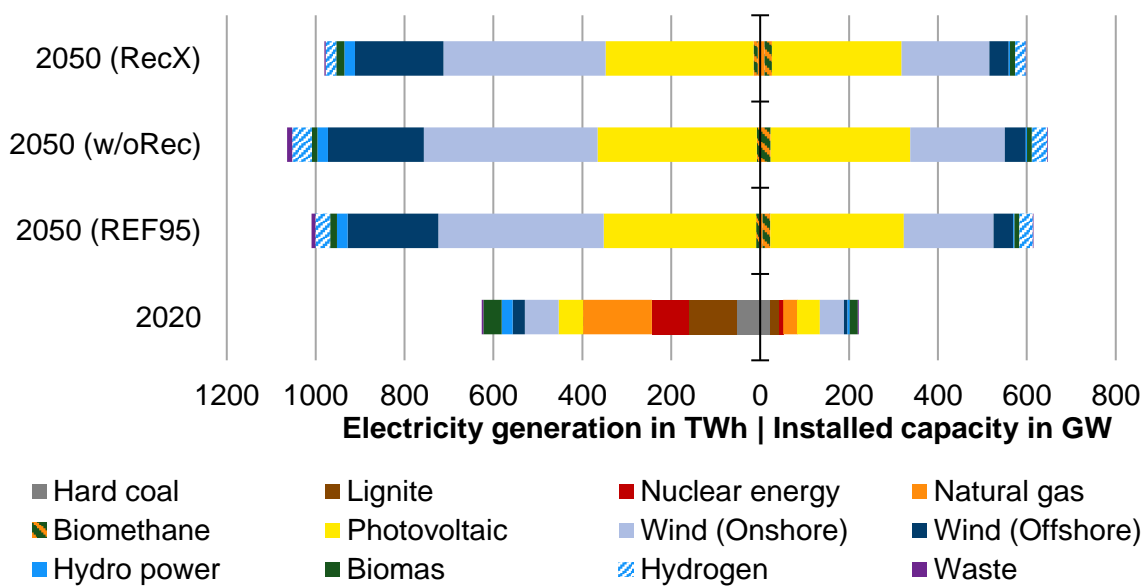

**Figure 3-3 Installed capacities and electricity generation across all scenarios in 2050 compared to 2020-values**

First, it can be stated that a significant expansion of renewable electricity generation is necessary in the reference scenario until the year 2050. More than 575 TWh of electricity will be generated by wind turbines and more than 340 TWh by photovoltaics in 2050. This requires 300 GW of photovoltaics, more than 200 GW of wind-onshore and 44 GW of wind-offshore to be installed. $CO_2$-intensive power generation from fossil fuels will be completely phased out by 2050. Increased recycling can save nearly 20 GW of installed capacity. Without recycling, more than 30 GW of additional capacity is needed. In conclusion, the impact on electricity demand and generation is relatively small compared to the Reference Scenario.

The situation is different when analyzing future hydrogen demand. In the reference scenario, approx. 400 TWh of hydrogen will be required by 2050, which will be used primarily in the industrial and transport sectors. With an increase in recycling rates, up to 125 TWh of hydrogen can be saved in 2050 (see Figure 3-4). On the other hand, a transformation of the energy system without recycling requires more than 170 TWh additional hydrogen in 2050. These results illustrate the great importance of recycling with its impact on the overall energy system, and especially with regard to hydrogen demand.

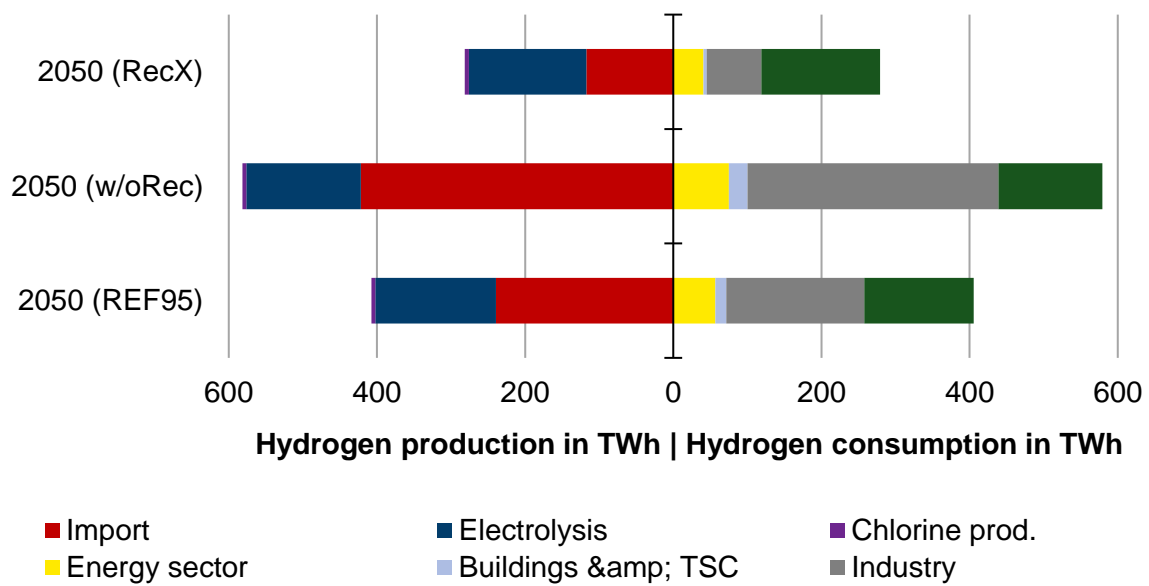

**Figure 3-4 Hydrogen production and consumption across all scenarios in 2050**

In addition, it must be noted here that if recycling measures are excluded, the increased demand for hydrogen will mainly be met by increased hydrogen imports. If hydrogen import quantities of this magnitude are not available for Germany, this quantity must be generated domestically with electrolysers, which means that different recycling strategies can ultimately have a major influence on the required electricity generation capacity. If future biomass imports are available, these could be exploited for example in the generation of industrial process heat or in the production of valuable chemical products (e.g., gasification).

### 3.2.2 Impact of recycling on the industrial sector

The following is a comparative scenario analysis for the industrial sector. The final energy demand is shown in Figure 3-5. The *w/oRec* and *RecX* scenarios illustrate the impact of recycling as an energy efficiency measure. Without recycling, industrial final energy consumption increases by about 300 TWh. On the other hand, an increase in recycling rates leads to a reduction in final energy demand of about 200 TWh. The biggest change can be observed in the demand for hydrogen. In scenario *w/oRec* the industry sector requires more than 340 TWh hydrogen in 2050 to compensate for the higher share of primary resource production. In comparison to the reference scenario an additional 65 TWh of natural gas are required for non-energetic demand.

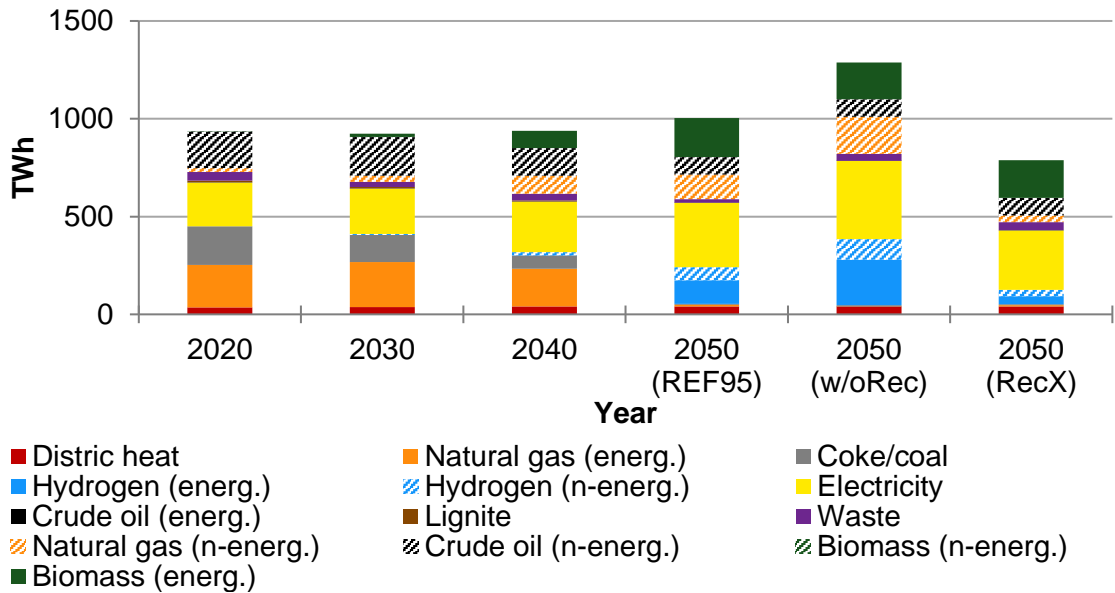

**Figure 3-5 Development of industrial final energy demand in reference scenario, with 2050-values across all scenarios**

Being able to exploit the theoretical maximum recycling rates decreases the hydrogen and natural gas demand in 2050 to half compared to *REF95*. The hydrogen consumption separated by industry in 2050 is illustrated in Figure 3-6.

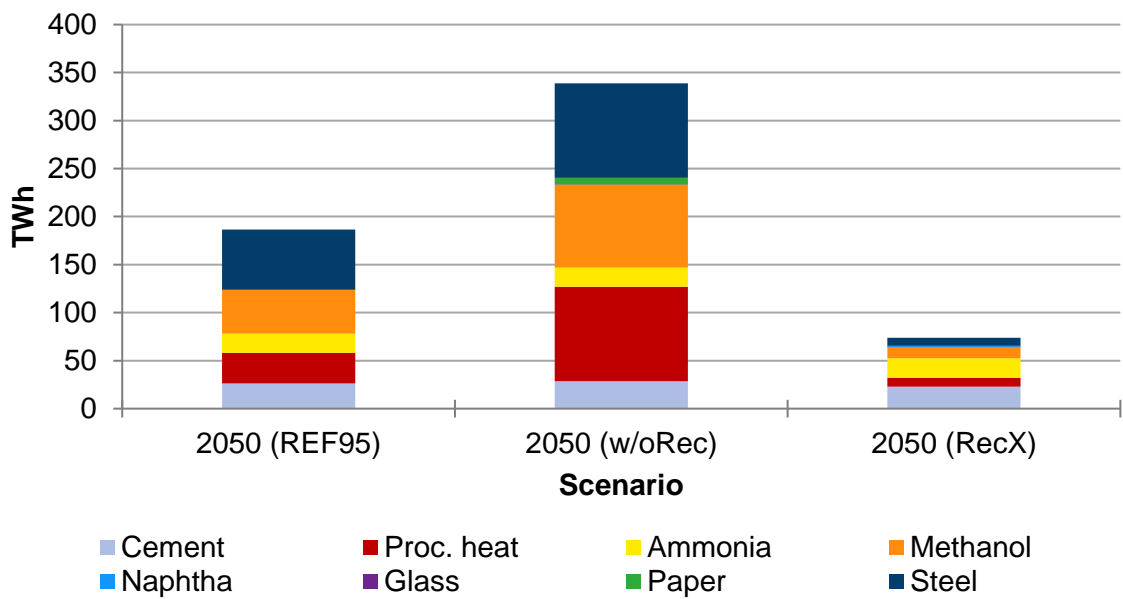

**Figure 3-6 Development of industrial hydrogen demand in reference scenario, with 2050-values across all scenarios**

Without recycling, more hydrogen is needed both for steel production and for the production of methanol. In addition, due to the lack of recycling measures, more process heat is required, which is provided by hydrogen furnaces. Increased recycling can reduce the demand for hydrogen in these industries and thus reduce the amount of hydrogen required for the entire energy system.

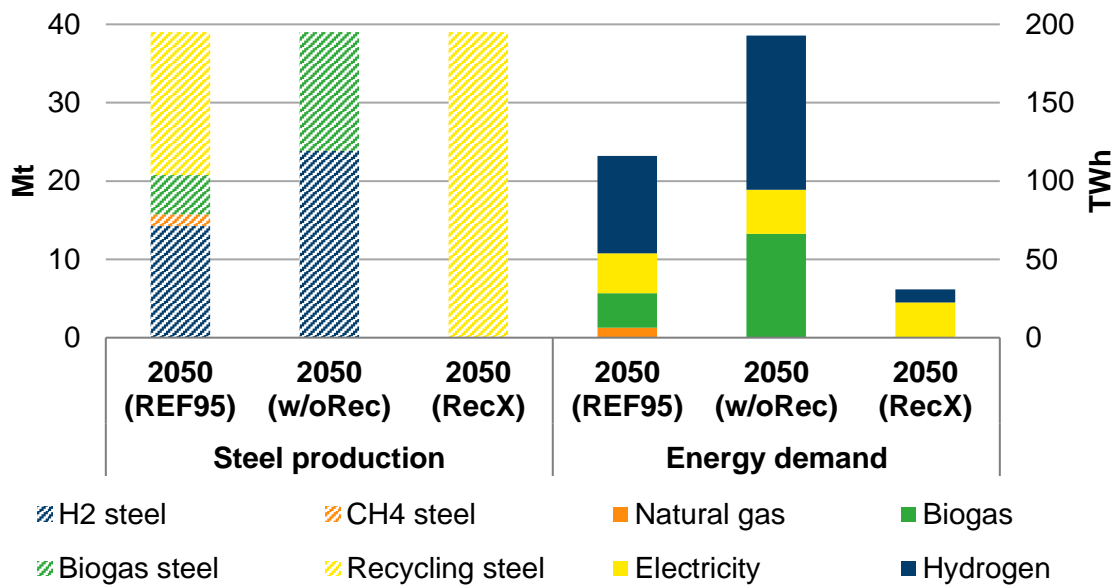

**Figure 3-7 Steel production and corresponding energy demand across all scenarios in 2050**

Figure 3-7 shows that steel recycling has a direct influence on the energy carrier mix in steel production. The more steel scrap that can be recycled, the less hydrogen and biogas are needed for direct reduction. Without recycling, the energy demand increases by about 80 TWh. It can be seen that the supply of hydrogen in the energy system is of significant importance for the success of the transformation of the energy system and that this can lead to competition in the demand for hydrogen between the chemical industry and steel production.

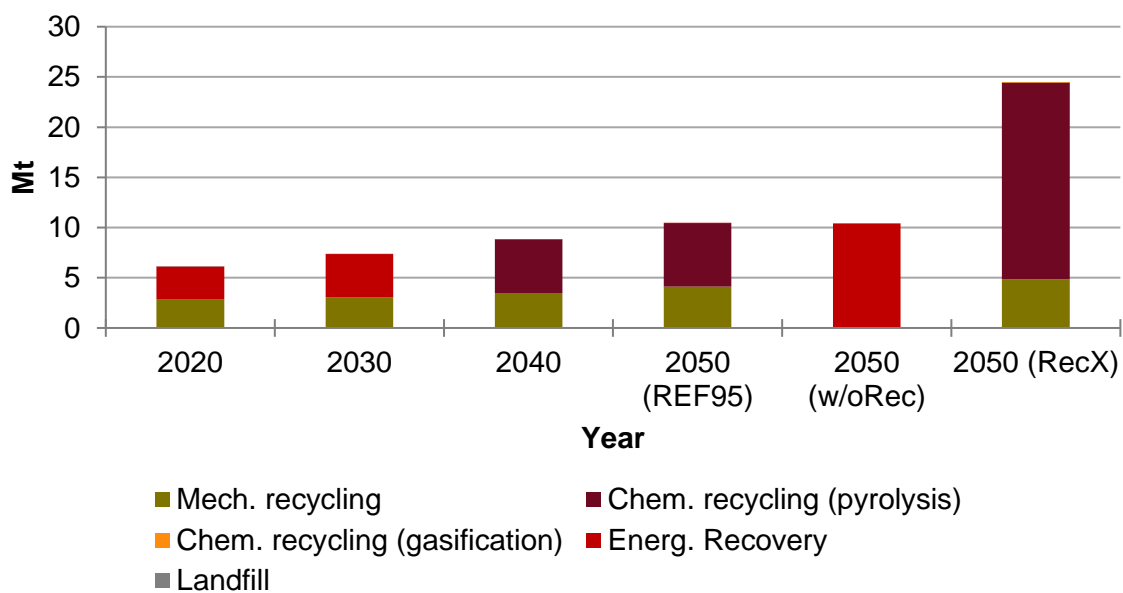

**Figure 3-8 Development of plastics waste processing in reference scenario, with 2050-values across all scenarios**

Currently, about 6 Mt of plastic and plastic waste is sent for recycling in Germany [26]. More than half of this is burnt with energy recovery, and the rest is mechanically recycled. The development of plastic waste recycling up to the year 2050 is shown in Figure 3-8. In the reference scenario, no more plastic waste will be burnt in 2050; instead, it will be exclusively mechanically or, for the most part, chemically recycled and thus returned to the material cycle. If the energy system model prohibits recycling into the feedstock, all plastic waste in 2050 must be recycled for energy (*w/oRec*). The end product of chemical recycling (pyrolysis oil) is therefore no longer available to the energy system and other raw materials must be used (see Figure 3-9). If, on the other hand, the energy system is left to decide how much plastic waste can be recycled, chemical recycling plays a significant role (*RecX*). Nearly 20 Mt of plastic waste will be chemically recycled in 2050 and converted into pyrolysis oil, which can then be used again to produce highly refined chemicals.

The production of high-valuable chemicals compared to the reference scenario is shown in Figure 3-9. For defossilization of the chemical industry, chemical recycling of end-of-life plastics to produce green naphtha (pyrolysis oil) and the methanol-to-olefins route are crucial. More than 10 Mt of high-valuable chemicals need to be produced additionally in both scenarios via the corresponding green variants.

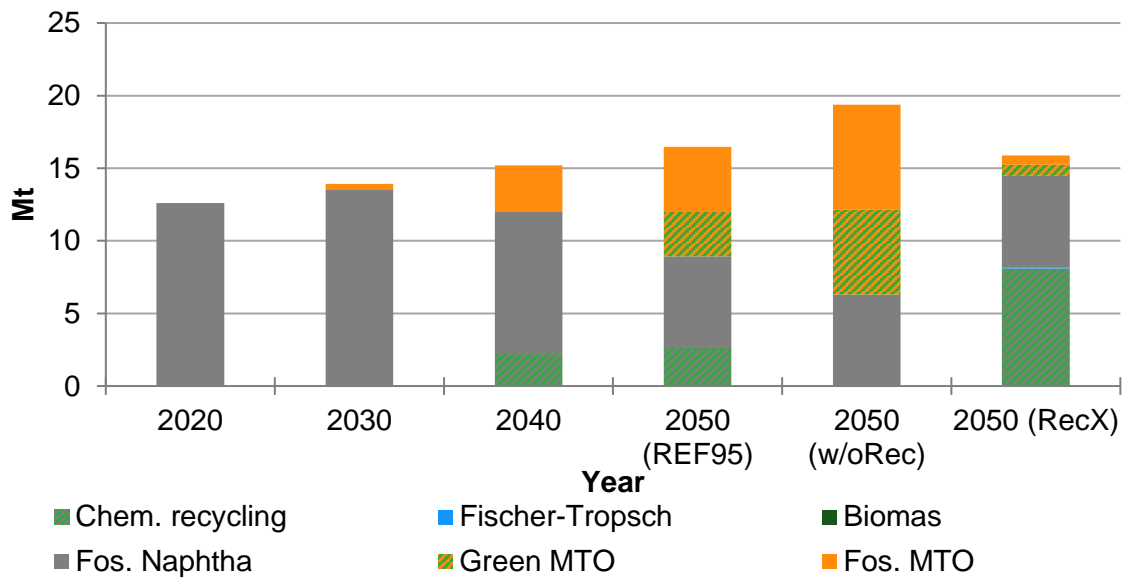

**Figure 3-9 Development of high-valuable chemicals production in reference scenario, with 2050-values across all scenarios**

In the *w/oRec* scenario, the omission of pyrolysis oil from chemical recycling must be substituted with fossil and green methanol. With increased recycling, an additional 5 Mt of high-valuable chemicals can be produced from pyrolysis oil. This will substitute both fossil methanol and green methanol. It can be concluded that chemical recycling of end-of-life plastics is more cost-effective than primary production via the methanol-to-olefins route.

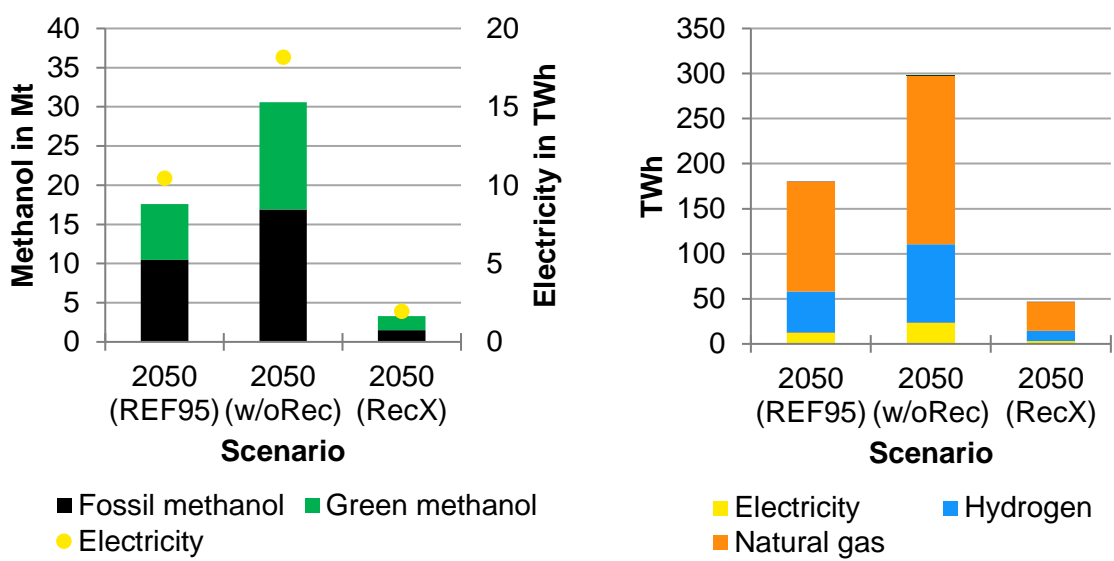

**Figure 3-10 Electricity and methanol demand for methanol-to-olefins route (left) and energy demand for methanol production (right) across all scenarios in 2050**

Nevertheless, the production of high-valuable chemicals via the methanol-to-olefins route will play a very decisive role in the future (see Figure 3-10). As an alternative to cracking naphtha, the methanol process is a significantly lower $CO_2$ variant. A ban on recycling means that approx. 30 Mt of methanol will be required in 2050 for the production of highly refined chemicals (*w/oRec*). This leads to the use of about 87 TWh of hydrogen and 187 TWh of natural gas in 2050. If higher recycling rates of end-of-life plastics are achieved, the total energy demand of methanol production decreases to about 50 TWh (*RecX*). This illustrates the great influence recycling has on the energy demand in industry, but also on the energy carrier mix.

### 3.2.3 Scrap costs sensitivity

Since an increase in recycling rates, as made possible in the *RecX* scenario, is not readily feasible (e.g., impurities in scrap, reduced quality of the steel products produced, availability of secondary raw materials), the following sensitivity analysis on the example of steel scrap costs was carried out. In this calculation, the costs in the *RecX* scenario for recycling processes and secondary raw materials are successively increased to determine the point up to which recycling is still worthwhile. For steel scrap, the model uses costs of 236 €/t steel scrap (comparable to today's costs [27]). The following Figure 3-11 shows the development of steel production in 2050 of the *RecX* scenario when the costs for steel scrap are gradually increased. It can be observed that a process change away from recycling to direct reduction with hydrogen only takes place from approx. 472 €/t steel scrap. From 570 €/t of steel scrap, steel recycling becomes too expensive for the energy system, so that more and more use is made of direct reduction with hydrogen, natural gas and biogas. From approx. 700 €/t steel scrap, steel recycling is no longer part of the model solution.

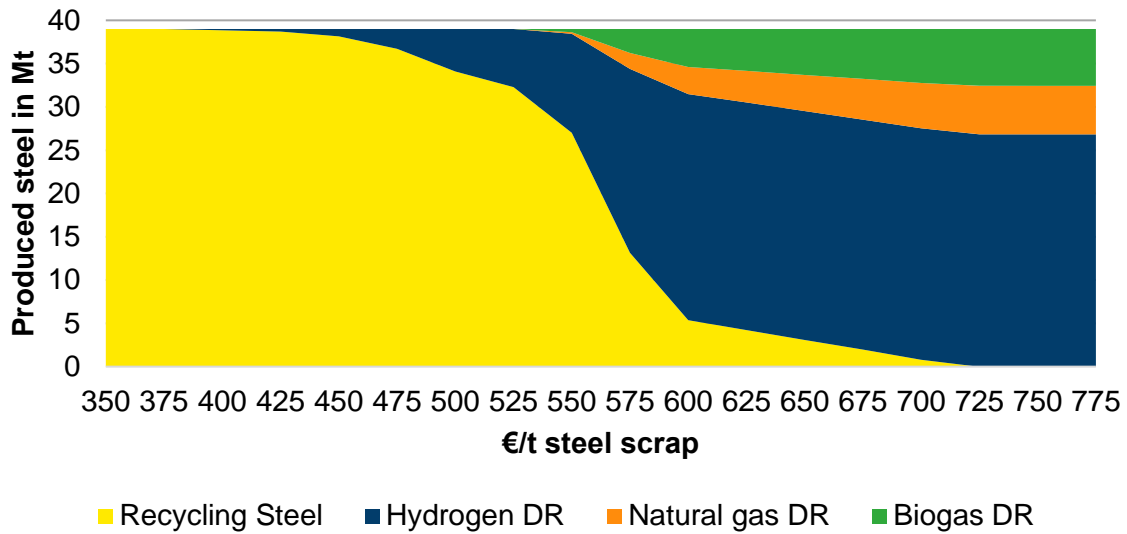

**Figure 3-11 Development of steel production as a function of steel scrap costs in 2050 in the *RecX* scenario (DR: Direct Reduction Process)**

The corresponding development of the energy demand of the steel industry can be seen in Figure 3-12. It should be noted here that the hydrogen demand, which is already apparent at low costs for steel scrap, is not used for direct hydrogen reduction, but to provide the heat required by the electric arc furnace.

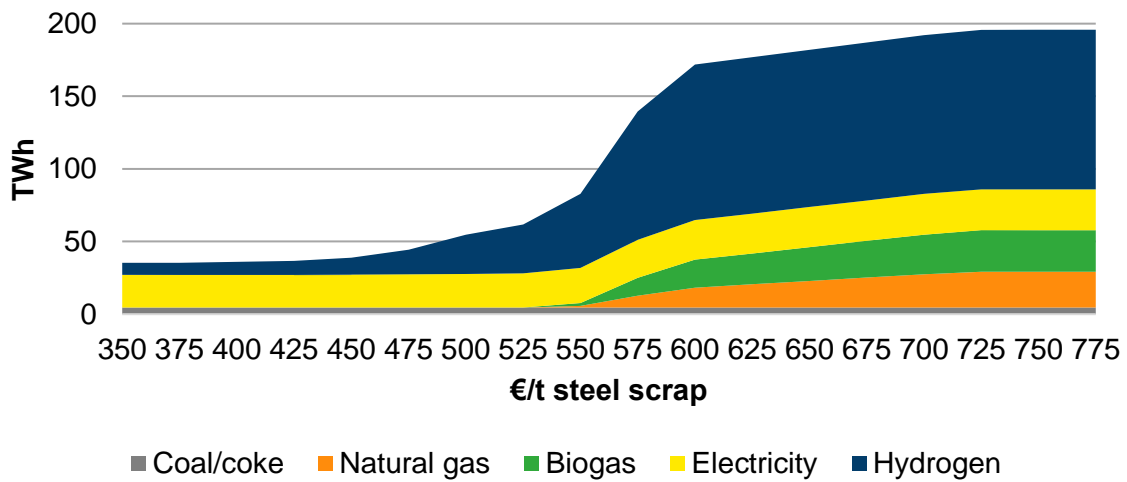

**Figure 3-12 Development of the energy demand in steel production as a function of steel scrap costs in 2050 in the *RecX* scenario**

Analogous to the development of steel production processes, a change in energy demand in the steel industry can only be observed when today's steel scrap costs double. A drastic increase in energy demand, essentially hydrogen, can even only be seen from a cost range of 580-590 €/t steel scrap.

It can be stated that steel recycling processes will only be replaced by other production processes (direct reduction) if steel scrap costs double compared to today's. Depending on the type of steel products, the steel scrap used must be pre-processed accordingly in order to maintain the required purity levels. The costs of reprocessing are difficult to quantify. Based on this sensitivity, the costs for pre-processing steel scrap can at least be as high as today's steel scrap costs and steel recycling would still be a cost-efficient $CO_2$ mitigation option.

### 3.2.4 Impact of recycling on the total system costs

The evaluation of the scenarios shows that recycling measures have a significant influence on the cumulative additional costs of the entire transformation (Table 3-2). Without recycling, the additional financial expenditure increases almost twofold. By contrast, maximum utilization of recycling rates has the potential to reduce the additional costs by more than a quarter compared with the reference scenario. The greatest effects can be observed in the industry and conversion sector, which includes investments in new processes and expansion of renewable energy technologies. Also, the import of renewable energy carriers (e.g., hydrogen) is highly dependent on recycling efforts and thus contributes to great cost differences between scenarios (Figure 3-13).

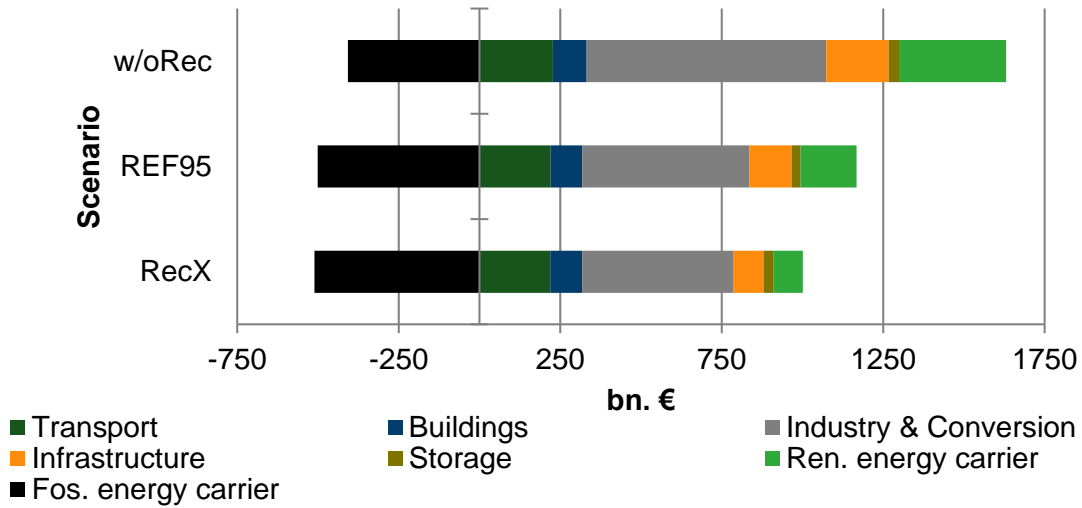

**Figure 3-13 Total cost of transformation cumulated from 2020 until 2050 across all scenarios**

It can be concluded that future efforts to increase recycling rates have a large leverage to reduce the monetary effort to achieve the climate targets (-26%). On the other hand, the cumulative additional costs nearly double (+84%) compared to the reference scenario if recycling measures are not available to the energy system. This leads to average $CO_2$ abatement costs over the entire transformation of 215 €/t$CO_2$ for the scenario without recycling (*w/oRec*) and 86 €/t$CO_2$ for the scenario with maximum recycling rates (*RecX*) (Table 3-2).

**Table 3-2 Cost data of recycling scenarios (comparison to reference scenario in parenthesis)**

|  | *w/oRec* | *RecX* |
|---|---|---|
| **Cost of transformation 2020-2050** | €1224 bn. (+557) | €490 bn. (-176) |
| **Accumulated $CO_2$-savings 2020-2050** | 5699 Mt$CO_2$ (±0) | 5699 Mt$CO_2$ (±0) |
| **Avg. spec. $CO_2$ abatement costs 2020-2050** | 215 €/t$CO_2$ (+98) | 86 €/t$CO_2$ (-31) |
| **Avg. spec. $CO_2$ abatement costs 2050** | 374 €/t$CO_2$ (+122) | 188 €/t$CO_2$ (-64) |
| **Marginal $CO_2$ abatement costs 2050** | 847 €/t$CO_2$ (+209) | 517 €/t$CO_2$ (-121) |

# 4 Conclusion

The aim of this work is to investigate the role of recycling measures as a greenhouse gas reduction strategy in the context of the overall energy system transformation. An integrated energy system model, which can be used to calculate cost-effective national transformation strategies for the German energy system is chosen as the basic model for the investigations in this paper. The metal and non-ferrous metal production, basic chemicals, paper industry, glass production, and processing of stones and earths were selected for the subsequent detailed modeling on the basis of their share in total emissions and final energy consumption. With a 66% share of industrial final energy consumption and more than 80% of industrial $CO_2$-emissions, these industries can be assumed to be representative of the German industrial sector. By analyzing these industries, conclusions can thus be drawn about the behavior of the industrial sector as a whole. The individual processes are wired in the model in such a way that they compete with each other, but also with technologies from other demand sectors. This integration allows for analysis of effects of whole-system $CO_2$ mitigation on the industrial sector as well as feedbacks of changes in the industrial sector on the whole system. Recycling options are implemented as an additional alternative to conventional processes that rely on primary feedstocks. As a result, they compete to meet the respective demand for goods in a cost-optimal manner while restricted by overall system $CO_2$-emission levels. While in many other studies recycling rates are only exogenously specified, they are now part of the optimization. This allows for the first time a cost-optimal evaluation of recycling in the context of the entire energy system. The amount of secondary raw materials available in the future is estimated using an approach from material flow modeling. A frequency distribution is used to calculate at which point in time a certain quantity of a material will be available to the energy system as waste in the future.

*Scenario analysis*

As a starting point for the various studies in this paper, a baseline scenario is created (*REF95*). In this scenario, a transformation strategy is calculated to achieve a $CO_2$ reduction of the energy system in 2050 by 95% (compared to 1990) in a cost-optimal way. The technologies and processes required for this serve as a benchmark for the scenarios on recycling strategies. Thus, the reference scenario is used to evaluate and rank the other scenarios. The recycling scenarios are divided into the scenario *w/oRec*, in which recycling measures are prohibited in the energy system model, and the scenario *RecX, in which* industrial goods can be produced entirely by recycling processes and secondary raw materials. The main findings of the scenario analysis are summarized below.

<u>Energy system without recycling (*w/oRec*)</u>

This scenario allows for the first time an estimation of the influence of recycling on the development of the energy system.

- Excluding the option to recycle results in an additional demand of 300 TWh in 2050 for primary energy demand, almost entirely due to higher energy demand in the industrial sector of 285 TWh.
- Overall, the industrial hydrogen demand is more than 150 TWh above that of *REF95* at about 350 TWh.

- A 95% reduction in $CO_2$ emissions in 2050 without recycling measures can only be achieved at considerable additional financial expense. In contrast to the reference case, the cumulative costs of the transformation increase by an additional 84%, or €557 billion.
- Furthermore, it can be stated that without recycling in today's energy system, additional costs of 13 billion €/a would arise.

Increased recycling rates (*RecX*)

For the first time, recycling rates are part of an energy system optimization and are not set exogenously. While only very rough estimates have been available so far, a detailed and consistent picture can now be drawn. The following results are to be underlined:

- By 2050, the possible recycling quotas will be utilized to the maximum. Only secondary raw materials will be used in the production of steel, aluminum, paper, glass and plastics.
- This results in a reduction of primary energy demand of 250 TWh and final energy demand in industry of 200 TWh in 2050 compared to the reference case.
- In 2050, steel production will require only 22 TWh of electricity for the electric arc furnace and no more hydrogen for direct reduction.
- In methanol production, both 90 TWh of natural gas and 34 TWh of hydrogen are eliminated because 14 Mt less methanol is processed into high-value chemicals in the methanol-to-olefins route. This is due to an increase in chemical recycling of end-of-life plastics, resulting in 8 Mt of primary high-value chemicals.
- Overall, in the *RecX* scenario, the additional costs of the transformation decrease by 26% (€176 billion) compared to *REF95*. Efforts to increase recycling rates have great potential to reduce the additional financial costs of the energy transition.

Recycling is of great importance for the energy system as a measure to mitigate transformation costs and should receive greater attention as an element of a greenhouse gas reduction strategy. The chosen model approach and the model philosophy are suitable to make statements about recycling in the context of GHG mitigation strategies. Recycling is a cost-effective measure in the context of greenhouse gas reduction strategies. The same holds true for decision makers, who need to pay more attention to recycling. Thus, future energy scenarios need to consider recycling and provide a more holistic view of the energy system design to enable more robust decision making.

This analysis is a first step to estimate the value of recycling within the context of energy system modeling. For future research it is of value to additionally analyze $CO_2$ footprints of imported materials. This would require broader system boundaries and global energy system models. Furthermore, the future availability of waste/scrap is very uncertain. Therefore, the combination of material flow models and energy system models should be strengthened.

# Appendix A  Basic assumptions

For the year 2050, the energy system model is provided with a limited amount of $CO_2$ that can continue to be emitted. However, exogenous emission targets for optimization must also be specified for the intermediate years. For the year 2030, based on the German government's Climate Protection Program 2030 [5] a 55% reduction in greenhouse gases compared to 1990 is specified. Since, apart from this target, no emission targets for further years are anchored in law [28] the emission targets for the intermediate years up to the 95% reduction in 2050 are interpolated linearly. The development of the exogenously set emission targets up to the year 2050 can be seen in Figure A-1. Sectoral targets are not specified, however, so that conclusions about the cost-effectiveness of certain measures across all sectors can also be drawn from the optimizations.

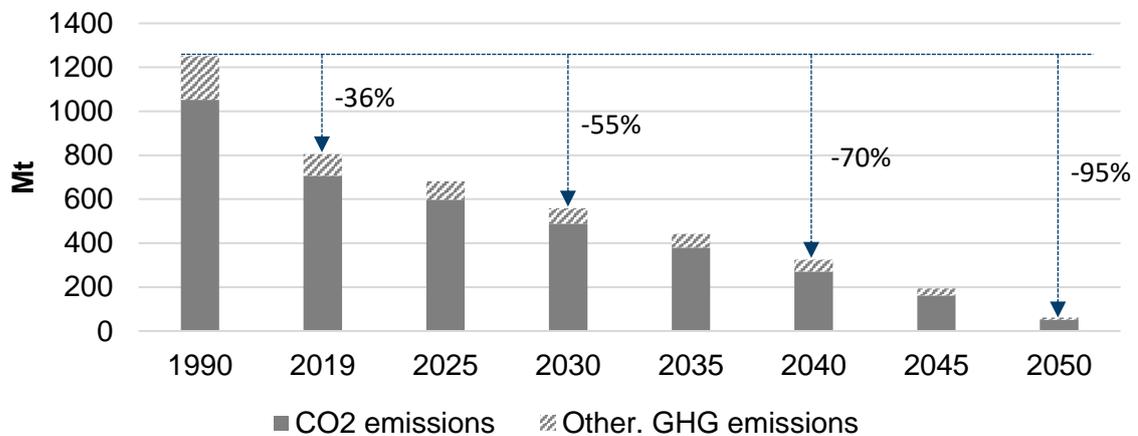

**Figure A-1 Development of exogenously set emission targets until 2050 in Mt**

A further basic assumption for the subsequent scenario calculations is the phase-out of coal-fired power generation by 2038 [25] and the phase-out of nuclear energy by 2022 [24]. Both are already enshrined in law and are therefore also used for the model calculations.

The historical and current values are available from the Federal Office of Economic [29–32] and the Federal Statistical Office [33] from the Federal Statistical Office. The future development of fuel costs corresponds to the scenarios from the World Energy Outlook [34].

In the final energy sectors, future demand is of particular importance. The development of the transport performance divided into the different means of transport is shown in Figure A-2. It is assumed that passenger transport performance will remain constant at around 1100 billion passenger kilometers until 2050, with the share of buses and passenger trains increasing slightly at the expense of private transport. Freight transport performance will increase by more than a third to 954 billion tonnage kilometers by 2050. Most of this will be done by trucks and freight trains. Today's figures are taken from the Federal Ministry of Transport in Figures and Digital Infrastructure [35]. The future progression is based on the study Climate Paths for Germany [12].

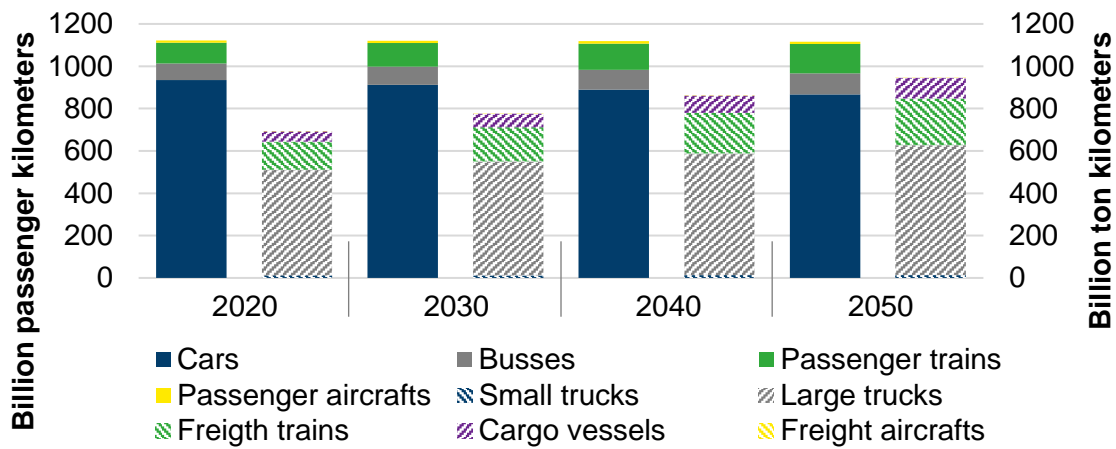

**Figure A-2 Development of exogenously set transport demand until 2050 (based on [12,35])**

In the building sector, both residential and non-residential buildings of the commercial, trade and service (CTS) sector as well as industry are combined in the model. The development of the demand for space heating and hot water divided among the different building types until 2050 is shown in Figure A-3. It is assumed that the demand of space heating and hot water in non-residential buildings will decrease. However, increasing demand in residential buildings will partially compensate for this decrease, resulting in a slight overall decrease in demand from about 800 TWh today to 760 TWh. The further subdivision of building types by building age classes are shown in Lopion et al. 2020 [17].

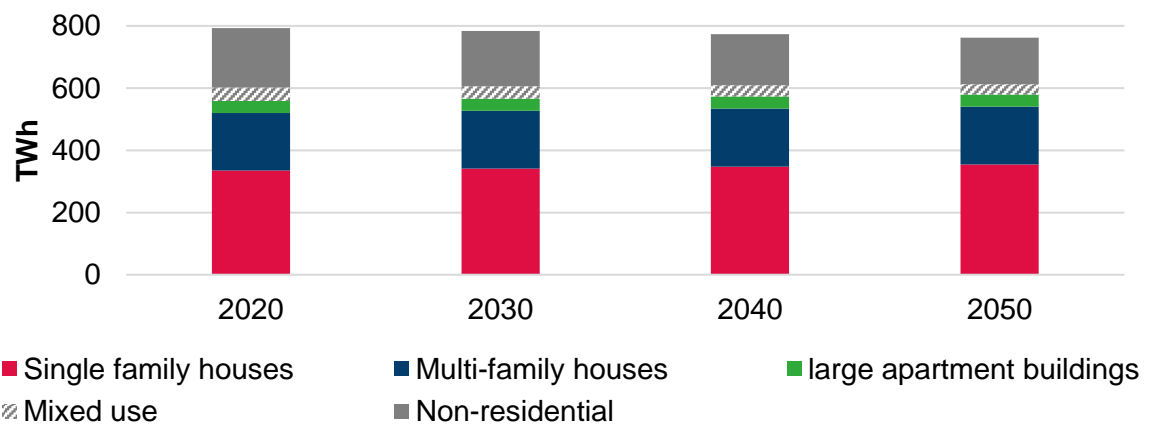

**Figure A-3 Evolution of exogenously set demand of space heating and hot water in residential and nonresidential buildings until 2050 (adapted from [12])**

Demand is also crucial for the industrial sector. The assumed demand for goods in the individual industrial sectors can be taken from Table A-1. For the remaining industry, which is not shown in detail, the development of the gross value added of the individual industrial sectors is decisive. The current production of goods is taken from the annual reports of the industrial associations or of the individual industries [36–41]. For all sectors, based on the Climate Paths for Germany study, it is assumed that energy demand will remain constant or increase [12]. It is assumed that industrial goods production will remain constant or increase up to the year 2050. Only steel production is expected to decline slightly from the current level of around 40 Mt of crude steel to 39 Mt.

**Table A-1 Development of demand for industrial goods until 2050 (based on [12])**

| in Mt | 2020 | 2030 | 2040 | 2050 |
|---|---:|---:|---:|---:|
| **Steel** | 40.0 | 40.0 | 39.0 | 39.0 |
| **Cement** | 34.2 | 35.1 | 35.9 | 36.9 |
| **Aluminum** | 1.3 | 1.3 | 1.3 | 1.3 |
| **Copper** | 0.7 | 0.7 | 0.7 | 0.7 |
| **Zinc** | 0.3 | 0.3 | 0.3 | 0.3 |
| **Urea** | 0.5 | 0.6 | 0.7 | 0.8 |
| **Ammonia** | 2.3 | 2.6 | 3.0 | 3.4 |
| **Methanol** | 1.1 | 1.3 | 1.5 | 1.7 |
| **Chlorine** | 3.9 | 4.5 | 5.1 | 5.7 |
| **Plastics** | 14.4 | 16.0 | 17.7 | 19.4 |
| **Glass** | 7.2 | 7.6 | 8.0 | 8.4 |
| **Paper** | 22.7 | 23.6 | 24.5 | 25.5 |
| **Copper** | 0.7 | 0.7 | 0.7 | 0.7 |
| **Zinc** | 0.3 | 0.3 | 0.3 | 0.3 |

## Acknowledgement

This work was supported by the Federal Ministry for Economic Affairs and Energy of Germany in the project METIS (project number 03ET4064A).

# The Value of Recycling for Low-Carbon Energy Systems - a Case Study of Germany's Energy Transition

## Supplementary Material (SM)

Felix Kullmann[1,*], Peter Markewitz[1], Leander Kotzur[1], Detlef Stolten[1,2]

[1] Institute of Energy and Climate Research, Techno-economic Systems Analysis (IEK-3) Forschungszentrum Jülich GmbH, Wilhelm-Johnen-Str., 52428 Jülich, Germany

[2] Chair of fuel cells, RWTH Aachen University, c/o Institute of Energy and Climate Research, Techno-economic Systems Analysis (IEK-3) Forschungszentrum Jülich GmbH, Wilhelm-Johnen-Str., 52428 Jülich, Germany

[*] Corresponding author. E-mail address: f.kullmann@fz-juelich.de (F. Kullmann).

## SM A: Model description

The energy system model NESTOR includes all system-relevant technologies of the German energy system, which are described by nodes and edges and built up as a large interacting network. The edges represent energy and material flows between the individual nodes (e.g., electricity and heat flows, iron ore flow). The nodes describe all energy and material sources, conversion technologies, energy and material storage, and energy and material sinks. Energy and material sources are components that bring a flow of energy or material into the system from outside (e.g., electricity import). All components that convert one type of energy or material flow into another belong to the class of conversion technologies (e.g., coal-fired power plant converts coal into electrical and thermal energy). Components that bring a flow of energy or material out of the system are called energy or material sinks (e.g., electricity export). In addition, there are components that can store energy and materials over multiple time steps (e.g., battery). Drivers for the interactions in this network are energy consumption relevant demands (e.g., demand for a certain passenger transport service, demand for industrial goods or demand for a certain space heating in residential buildings), which are exogenously given and are not part of the optimization. All previously listed components are related to each other through the network and fall under the influence of certain conditions [1].

The objective function of this optimization problem conditions the minimization of the annual system costs. Summed over all components, the annual fixed system costs are the product of installed capacity and average investment costs. The variable operating costs of the energy and material flows at a certain point in time can be summed up over all edges and time steps to the annual variable system costs. Both parts together result in the annual system costs, which are to be minimized [2]. The so-called penny-switching effect occurs in linear optimization problems. In this case, one technology is preferred over the other, even though it has only slightly lower costs. Although this behavior follows the logic of strict cost optimization, the results often seem implausible compared to the real market. Furthermore, the smallest technology cost changes can have a significant impact on the entire system design, so that the robustness of the obtained result is sometimes no longer given. To obtain more robust model results and avoid the penny-switching effect, Lopion et al. [2] added a quadratic part to the objective function. Here, the fixed annual cost of a technology no longer depends solely on the average investment cost, but on a range of investment costs. This model extension attempts to represent a system behavior that is more consistent with reality. The derivation and validation can be found in Lopion et al. [2]. In addition to the actual objective function of the optimization problem, the network composition is subject to additional constraints. For the generation of greenhouse gas reduction scenarios, an important constraint is the limitation of $CO_2$-emissions of the overall system, which can be used to specify greenhouse gas reduction targets for specific years (e.g., climate protection plan of the German government). Other constraints relate to the mathematical formulation of energy and material conservation and can be found in Lopion [3] for more information.

The energy system model is set up with an hourly resolution. This means that for each hour of an optimization year, generation, storage, conversion, and demand must meet the conditions described above. Since the model includes about 1000 technologies and many different time series, the resulting complexity of the optimization problem has a large impact on the computational effort. One way to reduce the computation time is to aggregate time series. The methodological concept from Kotzur et al. [4,5] and Hoffmann et al. [6,7] is applied for this purpose. The principle of time series aggregation is based on the fact that similar characteristics of a time series are grouped together and represented by less data points. Accordingly, the aggregated time series comprises fewer data entries than the original time series. Thus, only one representative period has to be optimized in the model calculation and not the larger original data series.

The NESTOR model has no real spatial resolution, and is therefore a so-called single-node model. This means that both energy production and energy demand are located at the same virtual point [8]. Things get complex when it comes to mapping energy infrastructures and renewable electricity generation. In order to deal with the problem of different feed-in time series of e.g. wind energy plants in different regions of Germany, pseudo-regions were designed [3]. The principle is to implement spatial characteristics only for energy infrastructures and fluctuating renewable power generation, all other system components remain untouched to keep the computational effort manageable [9]. Following Sanchis et al. [10] the German mainland is divided into seven regions for implementing spatial differences in onshore wind, rooftop photovoltaics, and open field photovoltaics. Two regions (North Sea and Baltic Sea) are created to map offshore wind energy installations. A detailed description of the approach can be found in Lopion et al. [3].

In addition to optimizing the German energy system for a single year, a myopic transformation path implemented in the NESTOR model, which allows for an analysis of the years between the current and the future optimized energy system in 2050. In the back casting approach, the energy system is first optimized for the target year 2050. This cost-optimized design of the future German energy system in 2050 is the defined target for the calculation of the transformation path starting from today's energy system. In 5-year steps, the years of the transformation path from 2020 to 2050 are subsequently optimized. The model expands individual technologies in certain expansion corridors in such a way that the target values are ultimately met in 2050 and, at the same time, the $CO_2$-emissions limits of the overall system are met in the intermediate years. The expansion corridor is largely determined by the development of existing plants and market penetration effects. A detailed description of the myopic transformation path analysis can be found in Lopion et al. [3].

# SM B: Industrial process modeling

**Metal Production**

The energy and mass balances, as well as the techno-economic parameters of pig iron and steel production can be found in Table SM B-1.

**Table SM B-1 Investment costs and specific energy requirements of selected steelmaking processes**

|  | Investment costs [a] | Specific energy demand [b] | | | | | Reference value |
|---|---|---|---|---|---|---|---|
|  | 2020 (2050) | Power | Coal/coke | Natural gas | Hydrogen | Process heat |  |
|  | €/t | kWh/kg | kWh/kg | kWh/kg | kWh/kg | kWh/kg |  |
| Blast furnace | 365 (365) | 0,062 | 5,057 | 0,072 |  |  | Pig iron |
| Oxygen converter | 128 (128) | 0,018 |  | 0,108 |  |  | Crude steel |
| Electric arc furnace | 184 (184) | 0,576 |  |  |  | 0,215 | Crude steel |
| $H_2$-Direct-reduction | 220 (220) | 0,127 |  |  | 1,808 | 1,516 | Iron sponge |

[a] adapted from [11–13]
[b] adapted from [12–20]

The investment costs available in the literature always refer to the new construction of a plant. For the fixed operating costs, 5% of the investment costs are assumed. The variable operating costs of the plants cannot be given as a lump sum, as the electricity and fuel costs are derived endogenously from the model. In the conventional blast furnace, approx. 1.42 kg $CO_2$ are produced per ton of crude steel produced. If electricity, process heat and hydrogen are generated from renewable sources in the energy system, the electric arc furnace and hydrogen direct reduction can be considered as almost $CO_2$-neutral. The process heat required for the electric arc furnace and the hydrogen direct reduction can be provided by both natural gas and hydrogen combustion. Natural gas combustion produces specific $CO_2$-emissions, which must be taken into account in the model.

**Non-ferrous metals**

The energy and mass balances, as well as the techno-economic parameters of the individual processes in the non-ferrous metals industry can be found in Table SM B-2. The processes in the nonferrous metals industry are already characterized by a high degree of electrification, which means that the $CO_2$ emissions generated during electricity generation are decisive for subsequent process route changes. The provision of the additional process heat required is not tied to a specific energy source in the model but is part of the optimization.

**Table SM B-2 Investment costs and specific energy requirements of selected processes in the non-ferrous metals industry**

|  | Invest-costs [a] | Specific energy demand [b] | | | | Reference value |
|---|---|---|---|---|---|---|
|  | 2020 (2050) | Power | Coal/Coke | Natural gas | Process heat |  |
|  | €/t | kWh/kg | kWh/kg | kWh/kg | kWh/kg |  |
| Hall-Heroult-Process | 5000 (5000) | 15,622 | 3,413 |  | 0,739 | Aluminum |
| Aluminum melting furnace | 500 (500) | 0,7811 |  |  |  | Aluminum |
| Hall-Heroult-Process (inert) | 5500 (5500) | 19,528 |  |  |  | Aluminum |
| Primary copper | 2873 (2873) | 1,436 | 0,101 | 1,016 |  | Copper |
| Kayser Recycling System | 1810 (1810) | 0,858 | 0,120 | 0,556 |  | Copper |
| Primary zinc | 1520 (1520 | 3,989 |  |  | 0,034 | Zinc |
| Rolling oxide process | 2000 (2000) | 4,078 | 0,15 | 0,051 |  | Zinc |
| Zinc smelting furnace | 315 (315) | 0,08 |  | 0,694 | 1,34 | Zinc |

[a] based on [21,22] and expert interview with *Trimet* Aluminium SE
[b] based on [23–31]

**Glass and fiber optics**

The basic material for glass production is a mixture of quartz sand, limestone and soda ash. This is melted at approx. 1500°C in melting tanks. The molten glass is then further processed, shaped and stress-relieved. The melting of the batch has both the highest energy requirement and the highest $CO_2$-emissions. Specifically, 0.75 kg of $CO_2$ is emitted per kg of glass ready for sale (estimate based on [32]). The main energy sources for the furnaces today are natural gas, fuel oil and electricity. In principle, however, other fuels can also be used for the provision of process heat [32]. Almost 3 Mt of waste glass was collected in 2019, of which about 85% could be recycled. This means that glass recycling has one of the highest recovery rates in the world. Waste glass in the form of cullet can be added to the raw batch in certain proportions to reduce the energy input for melting. Approximately 3% of the energy input can be saved for every 10% of cullet used in the batch [33]. In addition, the cullet input requires less primary raw materials, which are responsible for process emissions during melting.

A list of the melting processes considered and their techno-economic assumptions can be found in Table SM B-3. The furnaces for glass production can be divided based on the main energy source used. Today, natural gas, fuel oil and electricity are used. For the natural gas and fuel oil-fired furnaces, there is also the option of supplementary electric heating. In addition, in the case of natural gas-fired furnaces, oxygen can be used for firing instead of air (oxyfuel process).

**Table SM B-3 Investment costs and specific energy requirements of selected glass production processes**

|  | Invest-costs [a] | Specific energy demand [b] | | | | Reference value |
|---|---|---|---|---|---|---|
|  | 2020 (2050) | Power | Natural gas | Heating oil | Hydrogen |  |
|  | €/t | kWh/kg | kWh/kg | kWh/kg | kWh/kg |  |
| Furnace flat glass | 195 (195) | 0,10 - 1,10 | 0,36 - 1,90 | 1,49 - 1,77 | 0,20 - 2,09 | Flat glass |
| Furnace glass fiber | 195 (195) | 0,09 - 1,00 | 0,37 - 2,10 |  | 0,20 - 2,31 | Fiberglass |
| Furnace hollow glass | 195 (195) | 0,09 - 1,00 | 0,26 - 1,40 | 1,07 - 1,30 | 0,14 - 1,54 | Hollow glass |
| Furnace special glass | 195 (195) | 0,12 - 1,25 | 0,66 - 3,29 |  | 0,36 - 3,85 | Special glass |

[a] adapted from [22,34,35]
[b] Adapted from [36–44]

In the future, three more furnaces will be available to the model with a superboost (80% electricity, 20% natural gas), a hydrogen (100% hydrogen) and a gas-mix furnace (80% natural gas, 10% hydrogen, 10% electricity). Since the specific energy carrier input varies for each furnace and depending on the type of glass to be produced, average values across all furnaces are given.

**Pulp and paper**

Cellulose fibers obtained from wood in industrial processes are used as raw material for the production of paper, although theoretically any raw material containing cellulose would be suitable. Depending on how these fibers are prepared, they are referred to either as mechanical pulp, if wood has been mechanically defibered, or as chemical pulp, if wood has been chemically pulped. The use of mechanical pulp gives the end product greater stability, as lignins remain in the wood pulp during mechanical defibration, and is therefore mainly used for cardboard. Lignin is removed during chemical pulping, as this leads to yellowing of the paper when exposed to light and would be unsuitable for higher quality paper products [45]. The prepared fibers (paper pulp) are pressed into paper webs in a paper machine. During this process, the paper pulp, which initially has a water content of more than 90%, must be dried and dewatered. The greatest energy requirement in paper production is for pulping and subsequent pressing and drying [46]. Instead of wood as the primary raw material, fibers from recovered paper are used to a large extent today (for packaging paper > 90%). Due to high quality collection and collection systems of recovered paper, in 2019 so much pulp from recovered paper was used that secondary raw materials accounted for 78% of the raw material mix used [47]. Also due to these recycling rates, the specific energy consumption is on average 2.6 MWh/t paper [48]. A list of the paper production processes and their techno-economic assumptions can be found in Table SM B-4.

**Table SM B-4 Investment costs and specific energy requirements of selected papermaking processes**

|  | Invest-costs [a] | Specific energy demand [b] | | Reference value |
|---|---|---|---|---|
|  | 2020 (2050) €/t | Power kWh/kg | Process heat kWh/kg |  |
| Paper Machine Graphic Paper | 1300 (1300) | 0,458 | 1,019 | Graphic paper |
| Paper Machine Packaging Paper | 1058 (1058) | 0,144 | 1,019 | Packing paper |
| Paper machine Sanitary paper | 1705 (1705) | 0,971 | 1,019 | Sanitary paper |
| Paper Machine Special Paper | 1108 (1108) | 0,231 | 1,019 | Special paper |
| Mechanical Defibering | 300 (300) | 2,745 | -1,127 | Wood pulp |
| Chemical Defibering | 1355 (1355) | 0,702 | 5,859 | Fiber |

[a] adapted from [46,49]
[b] adapted from [50–52]

Due to the high share of biomass in the raw material input for paper production, biogenic residues and by-products account for a relatively high share of the process heat supply and additionally also of the on-site electricity generation. In addition, natural gas is used in industrial cogeneration plants to generate electricity and process heat. In the future, the use of other biogenic energy sources or hydrogen, for example, is conceivable and accordingly also part of the model optimization.

**Basic chemicals and refineries**

A small number of chemical and petrochemical products are responsible for more than 75% of $CO_2$ emissions in this sector [53]. The focus of the following analyses is therefore on the processes for the production of the chemical products chlorine, hydrogen, ammonia, urea, methanol, and the highly refined petrochemicals (HVC).

*Ammonia synthesis*

In Germany, approximately 2.5 Mt of ammonia was produced in 2019 [54], of which approx. 80% is used as fertilizer [55]. Ammonia is synthesized conventionally via the Haber-Bosch process. In a first step, hydrogen is synthesized via natural gas steam reforming. In this process, natural gas reacts with steam at approx. 800°C to form a synthesis gas consisting of carbon monoxide and hydrogen. By adding compressed air, both the nitrogen contained in the synthesis gas and carbon monoxide are oxidized to carbon dioxide, which can then be separated more easily as carbon monoxide. The purified synthesis gas consisting of nitrogen and hydrogen then reacts in a high-pressure tank at approx. 450°C to produce ammonia [56]. In total, this process requires approx. 5.8 MWh natural gas/t ammonia for hydrogen production and a further 1.8 MWh/t ammonia in the form of process heat and approx. 2 MWh electricity/t ammonia. The resulting process emissions amount to 1.2 t$CO_2$/ton of ammonia. [57]. Since these emissions are due to steam reforming, a future hydrogen synthesis without process- and energy-related $CO_2$-emissions (e.g. water electrolysis with green electricity) would be an emission-neutral alternative. A list of selected processes for the production of important basic chemicals and their techno-economic assumptions (based on [53,57,58]) can be found in Table SM B-5.

**Table SM B-5 Investment costs and specific energy requirements of selected production processes in the chemical industry**

|  | Invest-costs [a] | Specific energy demand [b] | | | | Reference value |
|---|---|---|---|---|---|---|
|  | 2020 (2050) | Power | Natural gas | Hydrogen | Process heat |  |
|  | €/t | kWh/kg | kWh/kg | kWh/kg | kWh/kg |  |
| Chlor-alkali electrolysis | 404 (404) | 2,35 |  |  | 0,3 | Chlorine |
| Haber-Bosch process | 670 (670) | 2,07 | 5,83 |  | 1,83 | Ammonia |
| Haber-Bosch process ($H_2$) | 500 (500) | 1,72 |  | 5,93 | 1,83 | Ammonia |
|  |  | Power | Natural gas | Crude oil | Process heat |  |
|  |  | kWh/kg | kWh/kg | kWh/kg | kWh/kg |  |
| Methanol (steam reforming) | 400 (400) | 0,17 | 6,94 |  | 3,14 | Methanol |
| Methanol (partial oxidation) | 530 (530) | 0,18 |  | 9,22 |  | Methanol |
|  |  | Power | Hydrogen | Biomass | Process heat |  |
|  |  | kWh/kg | kWh/kg | kWh/kg | kWh/kg |  |
| Methanol ($H_2$) | 197 (197) | 1,5 | 6,33 |  |  | Methanol |
| Methanol (biomass) | 400 (400) | 0,17 |  | 10,08 | 3,33 | Methanol |

[a] based on [57]
[b] based on [53,58]

In ammonia synthesis, both natural gas (hydrogen production from steam reforming) and hydrogen provided directly externally (from within the energy system) can be used. In this process, natural gas and hydrogen are used raw materially, i.e. non-energetically. Additionally required process heat is also provided by natural gas in the conventional Haber-Bosch process, but can also be substituted by renewable energy sources in the future. In methanol synthesis, both externally provided hydrogen and biomass can be used instead of natural gas (hydrogen production from steam reforming) or heavy oil ((hydrogen production from partial oxidation). This energy carrier input is also a non-energy demand. However, since the feedstock use of externally supplied hydrogen lacks the carbon required for methanol synthesis (1.37 kg$CO_2$/kg of methanol), it must be additionally supplied. Other processes depicted in the energy system model in which $CO_2$ can be used feedstock-wise can be found in Table SM B-6. Naphtha is traditionally used to produce high-value chemicals.

**Table SM B-6 Investment costs and specific energy requirements of selected processes for the production of high-value chemicals (HVC) and synthetic refinery products**

| | Invest-costs [a] | Specific energy demand [b] | | | | Refer-ence value |
|---|---|---|---|---|---|---|
| | 2020 (2050) | Power | Methanol | Naphtha | Hydrogen | |
| | €/t | kWh/kg | kg/kg | kg/kg | kWh/kWh | |
| Methanol-to-Olefins | 268 (268) | 1,39 | 2,34 | | | HVC |
| Steamcracker | 1700 (1700) | 0,1 | | 1,22 | | HVC |
| Steamcracker (el. heating) | 250 (250) | 4,7 | | 1,22 | | HVC |
| Fisher-Tropsch-Synthesis | 788 (500) | kWh/kWh 3,5 | | | 1,32 | Diesel/ Gasoline/ Kero-sene |

[a] adapted from [59–61]
[b] based on [57,62,63]

In the future, however, it will also be possible to feedstock methanol and produce highly refined chemicals via methanol-to-olefins. Methanol produced by green hydrogen can contribute to the defossilization of petrochemicals. Other refinery products (diesel, gasoline, kerosene) can be produced using Fischer-Tropsch synthesis. The synthesis gas required for this consists of hydrogen and carbon monoxide, whereby the required carbon component can be provided from $CO_2$.

**Other industrial processes**

Industries and processes that are not covered by the industrial processes described above are mapped at a higher aggregation level. This is done using the methodology from Lopion et al. [3]. For this purpose, the current process heat demand of these industries, divided into three temperature levels, as well as the electricity demand is required (cf. Table SM B-7).

**Table SM B-7 Process heat and electricity demand in 2019 for aggregated mapped industrial processes in TWh (adapted from [3,64–66])**

| Branch | Process heat in TWh | | | Electricity in TWh |
|---|---|---|---|---|
| | < 100°C | 100-500°C | > 500°C | |
| Stone and earth industry | 0,03 | 0,05 | 2,34 | 1,72 |
| Nutrition & Tobacco | 18,11 | 22,25 | 0,00 | 19,08 |
| Other. chemical industry | 2,53 | 3,93 | 10,48 | 6,78 |
| Rubber & plastic goods | 1,74 | 6,59 | 0,00 | 13,81 |
| Metalworking | 4,11 | 3,25 | 6,00 | 14,64 |
| Mechanical Engineering | 2,92 | 2,18 | 4,12 | 10,89 |
| Vehicle construction | 5,11 | 3,88 | 7,15 | 16,06 |
| Other manufacturing | 18,18 | 6,07 | 7,27 | 21,14 |

No specific industrial goods are set as industrial demand, but rather the gross value added generated by the individual industry. Accordingly, a specific process heat and electricity demand per generated gross value added is stored in the model. The future energy demand of these industries is scaled up to the year 2050 based on forecasts of gross value added. The development of gross value added and thus the future energy demand of the sectors is taken from the study Climate Paths for Germany [67]. The process heat and the required electricity can be provided in the model by various technologies and energy carriers of the overall energy system. A detailed list can be found in Lopion et al. [3].

Unlike for the industrial processes shown in detail, no statements can be made about future process changes or process-specific energy efficiency gains due to economic decisions for the aggregated industries.

## SM C: Assumptions on future secondary raw material availability

### Steel

The estimation of future steel scrap quantities is based on the calculations from the material flow model in Pauliuk et al. [68] . This model is based on assumptions regarding the sectoral distribution of steel demand and the residence times stored in these stocks, which can be found in Table SM C-1. In doing so, Pauliuk et al. also estimate the historical steel production in Germany and are thus able to make statements on the steel supply in the four sub-divided sectors. Historical imports and exports of goods containing steel and steel scrap are also taken into account. As a simplification, the historical import and export rates are extrapolated to the year 2050.

**Table SM C-1 Assumptions on steel flows in Germany according to Pauliuk et al. [68]**

|  | Transportation | Mechanical Engineering | Construction | Other Products |
|---|---|---|---|---|
| **Sectoral breakdown** | 0,3 | 0,1 | 0,47 | 0,13 |
| **Average lifetime in years** | 13 | 20 | 50 | 10 |
| **Share of obsolete stock** | 0 | 0 | 0,1 | 0 |
| **Recovery rate** | 0,82 | 0,87 | 0,82 | 0,58 |

In addition to the sectoral breakdown of German steel demand and the average dwell time in these sectors, Pauliuk et al. also give an obsolete stock and recovery rate. The amount of steel that remains in a given sector and cannot be recovered is referred to as obsolete. It is assumed that approximately 10% of the steel flowing into the construction sector is no longer available for potential recovery. This takes into account the amount of steel that is installed in foundation walls or underground infrastructure and generally remains there [69]. Furthermore, a recovery rate is assumed for the sectors, which indicates what proportion can be recovered from the steel scrap flow of the respective sectors and what proportion is lost together with other scrap. This results in a theoretically maximum available quantity of steel scrap by the year 2050.

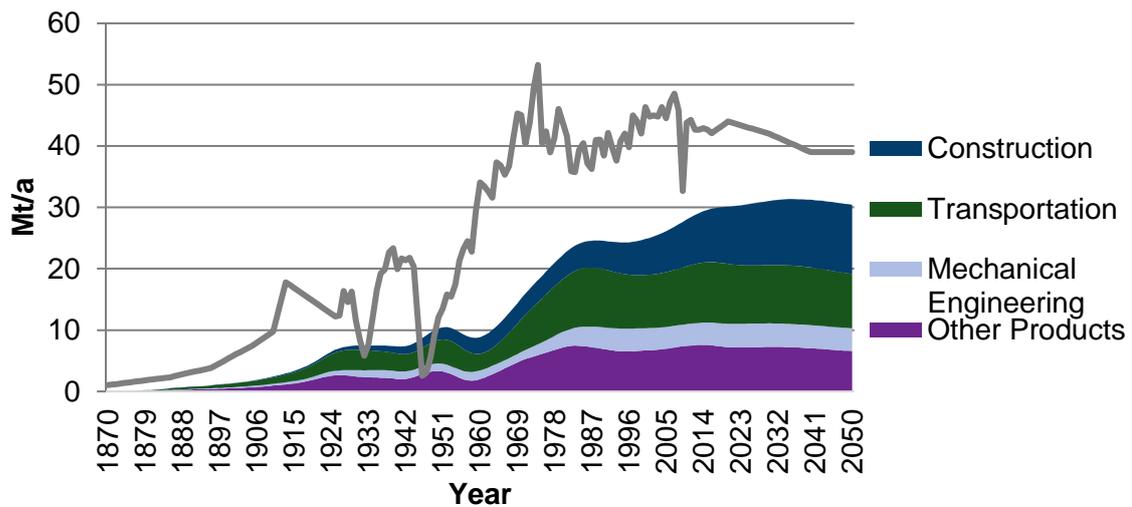

**Figure SM C-1 Available steel scrap quantity in Mt/a by 2050, adapted from Pauliuk et al. [68]**

Figure SM C-1 shows the historical steel production since 1870 and the amount of steel scrap generated annually from the four selected sectors: construction, transportation, mechanical engineering and other products. By 2050, approximately 30 Mt of steel scrap could be recovered from scrap production in Germany alone. It should be noted that this estimate does not provide any information on the economic viability of recovering steel scrap in the future.

**Non-ferrous metals**

A similar material flow analysis has been performed by Manfredi et al. [70] on behalf of the European Commission for aluminum and copper. The material flows and stocks for these non-ferrous metals within the EU until the year 2040 were estimated using the Global Aluminum Flow Model [71]. The assumptions underlying the model can be found in Table SM C-2.

**Table SM C-2 Assumptions on aluminum flows in Europe according to Manfredi et al. [70]**

|  | Transportation | Mechanical Engineering | Construction | Other Products |
|---|---|---|---|---|
| **Sectoral breakdown** | 0,29 | 0,09 | 0,26 | 0,36 |
| **Average length of stay in years** | 20 | 40 | 50 | 12 |
| **Recovery rate** | 0,94 | 0,64 | 0,86 | 0,66 |

It should be noted here that Manfredi et al. assume that approx. 12% of the aluminum produced in Europe is used as packaging material. Due to the short life of packaging, it is assumed in simplified terms that this quantity will already be produced again as aluminum scrap in Europe within one year. Due to the fact that material flow models rarely provide data up to the year 2050 and exclusively for Germany, an additional adjustment must also be made for the estimation by Manfredi et al. The results of the material flow analysis of Manfredi et al. are scaled down proportionally, taking into account the historical aluminum production in Germany, and extrapolated to the year 2050 under the assumption of constant production. The historical aluminum production for Germany can be found in the data of the statistical yearbooks of the Bureau of Mines of the United States of America [72]. Figure SM C-2 gives an example of the annually available aluminum scrap quantity for Germany until the year 2050, based on the calculations of Manfredi et al. [70].

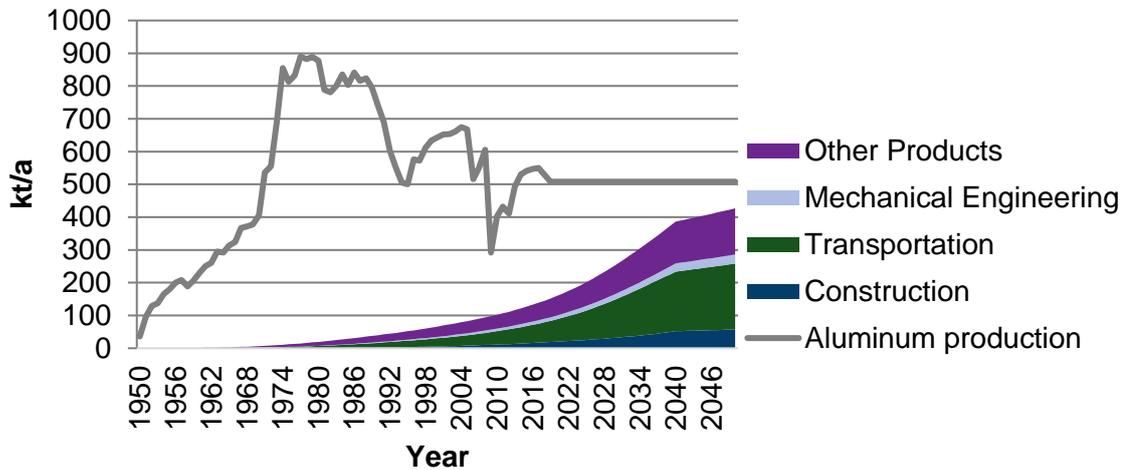

**Figure SM C-2 Available aluminum scrap quantity in kt/a until 2050, adapted from Manfredi et al. [70]**

These estimates only represent a theoretical maximum amount of steel scrap available, which do not allow any conclusions to be drawn regarding economic utilization.

**Plastics**

For the estimation of future plastic scrap volumes in Germany, this paper uses a study by Conversio Market & Strategy GmbH, as the author is not aware of any other holistic material flow analyses for future plastic flows in Germany. This study draws a material flow picture of plastic flows in Germany for the year 2019 [73]. Based on this, plastic waste volumes across five selected sectors result for the year 2019, which can be found in Table SM C-3.

**Table SM C-3 Plastic consumption and waste generation of selected sectors in kt in 2019, adapted from Conversio [73]**

| 2019 | Packing | Construction | Transportation | Electronics | Other Products |
|---|---|---|---|---|---|
| **Plastics consumption in kt** | 4369 | 3583 | 1509 | 881 | 3893 |
| **Waste volume in kt** | 3081 | 495 | 232 | 307 | 1119 |

Together with the assumptions on the development of future plastics production from Gerbert et al. [67] (+1.3%/a), this results in annual waste volumes up to the year 2050, which serve as input data for the energy system model. The development of the available waste volumes is shown in Figure SM C-3.

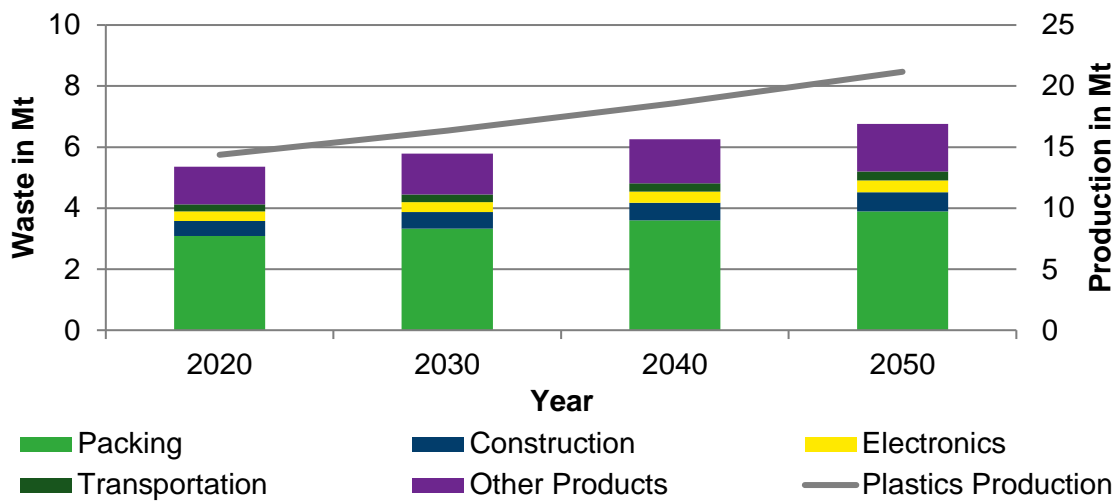

**Figure SM C-3 Available plastic waste quantities in Mt by 2050, adapted from [73][57][67]**

Geres et al. assume an annual increase of approx. 0.9%/year for the development of waste quantities. The resulting waste volumes are potentially available for recycling processes or energy recovery by the year 2050. Which share is used in the respective processes is part of the cost optimization of the energy system model.

**Glass**

For the estimation of future waste glass flows in Germany, this work is based on the historical values of the Federal Statistical Office [74]. These values are extrapolated taking into account an annual increase of 0.5%/a (according to the BDI study [67]) are extrapolated up to the year 2050. The proportion of the amount of used glass that could actually be recovered ranged from 82% (2008) to 89% (2014). Therefore, an average recovery rate of 85% is assumed for the transformation until 2050. The amount of waste glass deposited for the energy system model up to the year 2050 is shown together with the assumed glass production in Figure SM C-4.

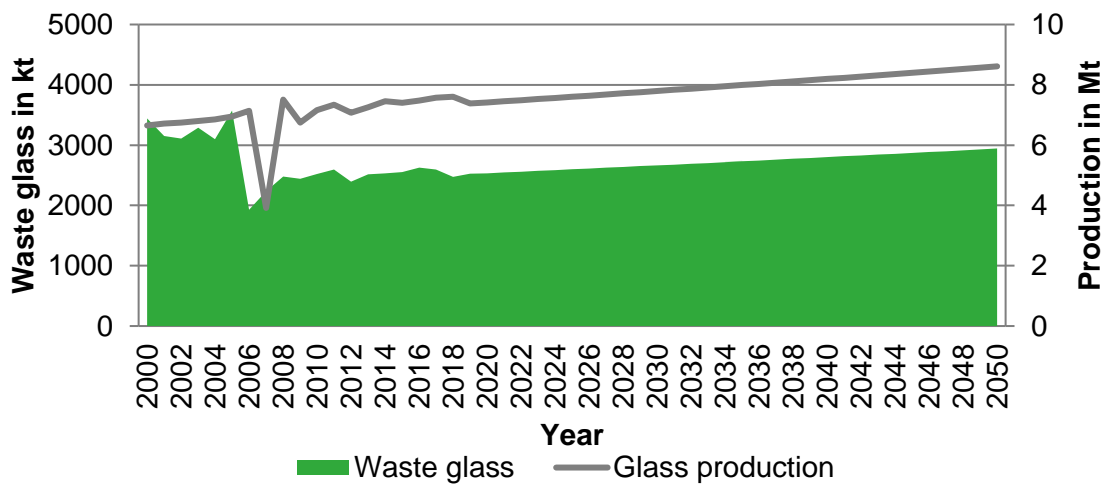

**Figure SM C-4 Waste glass volume in kt by 2050 according to own calculations based on Federal Statistical Office [74]**

In 2050, it is assumed that approx. 3000 kt of waste glass will be available for recycling. No detailed differentiation of glass qualities is made. These quantities of waste glass thus theoretically define maximum recycling potentials in the energy system model.

**Paper**

For paper products, due to their short residence time in the anthropogenic stock, it is assumed that the amount of recovered paper is directly available again for recycling processes in the energy system model. This results in a direct correlation of available recovered paper to annual paper production. The recovery rate for recovered paper over the last five years has averaged 75% [47]. This value is therefore also assumed for the transformation up to the year 2050. In 2050, taking into account an annual production increase of 0.4%/a (according to the BDI study [67]), the amount of paper produced will be 25.6 Mt, of which approx. 19.2 Mt of recovered paper can be used.

## Acknowledgement

This work was supported by the Federal Ministry for Economic Affairs and Energy of Germany in the project METIS (project number 03ET4064A).